\documentclass[a4paper,reqno, doublespace]{amsart}
\usepackage{eurosym}
\usepackage{amsmath, amssymb}
\usepackage{amsfonts}
\usepackage{graphicx}
\usepackage{color}
\usepackage[usenames,dvipsnames]{xcolor}
\usepackage[english]{babel}
\usepackage{bm}

\newtheorem{teo}{Theorem}[section]

\numberwithin{equation}{section}
\newtheorem{Lemma}{Lemma}[section]
\newtheorem{theorem}[Lemma]{Theorem}

\newtheorem{lemma}[Lemma]{Lemma}

\newtheorem{remark}[Lemma]{Remark}
\newtheorem{definition}[Lemma]{Definition}

\begin{document}
\title[NSC]{A boundary control problem for stochastic 2D-Navier-Stokes
equations}
\author{Nikolai Chemetov}
\address{Department of Computing and Mathematics, University of S{\~a}o
Paulo, 14040-901 Ribeir{\~a}o Preto - SP, Brazil}
\email{nvchemetov@gmail.com }
\author{Fernanda Cipriano}
\address{Center for Mathematics and Applications (NOVA Math) and Department
of Mathematics, New University of Lisbon, Portugal}
\email{cipriano@fct.unl.pt}
\date{}

\begin{abstract}
We study a stochastic velocity tracking problem for the 2D-Navier-Stokes
equations perturbed by a multiplicative Gaussian noise. From a physical
point of view, the control acts through a boundary injection/suction device
with uncertainty, modeled by stochastic non-homogeneous Navier-slip boundary
conditions. We show the existence and uniqueness of the solution to the
state equation, and prove the existence of an optimal solution to the
control problem.
\end{abstract}

\maketitle

%\tableofcontents



\section{Introduction}

\setcounter{equation}{0}

The goal of this article is to study an optimal boundary control problem for
stochastic viscous incompressible fluids, filling a bounded domain ${%
\mathcal{O}}\subset \mathbb{R}^{2}$, and governed by the Stochastic
Navier-Stokes equations with non-homogeneous Navier-slip boundary conditions%
\begin{equation}
\left\{ 
\begin{array}{ll}
\begin{array}{l}
d\mathbf{y}=(\nu \Delta \mathbf{y}-\left( \mathbf{y\cdot }\nabla \right) 
\mathbf{y}-\nabla \pi )\,dt+\mathbf{G}(t,\mathbf{y})\,d{\mathcal{W}}_{t}, \\ 
\\ 
\mbox{div}\,\mathbf{y}=0,%
\end{array}
& \mbox{in}{\ \mathcal{O}_{T}=(0,T)\times \mathcal{O}},\vspace{2mm} \\ 
\mathbf{y}\cdot \mathbf{n}=a,\;\quad \left[ 2D(\mathbf{y})\,\mathbf{n}%
+\alpha \mathbf{y}\right] \cdot {\bm{\tau }}=b\;\quad & \text{on}\ \Gamma
_{T}=(0,T)\times \Gamma ,\vspace{2mm} \\ 
\mathbf{y}(0,\mathbf{x})=\mathbf{y}_{0}(\mathbf{x}) & \mbox{in}\ {\mathcal{O}%
},%
\end{array}%
\right.  \label{NSy}
\end{equation}%
\vspace{1pt}where $\mathbf{y}=\mathbf{y}(t,\mathbf{x})$ is the 2D-velocity
random field, $\pi =\pi (t,\mathbf{x})$ is the pressure, $\nu >0$ is the
viscosity and $\mathbf{y}_{0}$ is the initial condition that verifies 
\begin{equation}
\mbox{div}\,\mathbf{y}_{0}=0\qquad \mbox{ in   }\ {\mathcal{O}}.
\label{ICNS}
\end{equation}%
Here $\ $%
\begin{equation*}
D(\mathbf{y})=\frac{1}{2}[\nabla \mathbf{y}+(\nabla \mathbf{y})^{T}]
\end{equation*}%
is the rate-of-strain tensor; $\mathbf{n}$ is the external unit normal to
the boundary $\Gamma \in C^{2}$ of the domain ${\mathcal{O}}$ and $\bm{\tau }
$ is the tangent unit vector to $\Gamma ,$ such that $(\mathbf{n},\bm{\tau }%
) $ forms a standard orientation in $\mathbb{R}^{2}.$\ The positive constant 
$\alpha $ \ is the so-called friction coefficient. The quantity $a$
corresponds to the inflow and outflow fluid through $\Gamma $, satisfying \
the compatibility condition%
\begin{equation}
\int_{\Gamma }a(t,\mathbf{x})\,\,d\mathbf{\gamma }=0\quad \quad 
\mbox{ for
any  }\;t\in \lbrack 0,T].  \label{eqC2}
\end{equation}%
This condition means that the quantity of inflow fluid should coincide with
the quantity of outflow fluid. The boundary functions $a$ \ and $b$ will be
considered as the control variables for the physical system \eqref{NSy}. The term $\mathbf{G}(t,\mathbf{y})\,{%
\mathcal{W}}_{t}$ is a multiplicative white noise.

The main goal of this paper is to control the solution of the system (\ref%
{NSy}) by the boundary condition $(a,b)$, which is a predictable stochastic
process belonging to the space $\mathcal{A}$ of admissible controls to be
defined in Section \ref{sec3}. The cost functional is given by 
\begin{equation}
\displaystyle J(a,b,\mathbf{y})=\mathbb{E}\int_{{\mathcal{O}}_{T}}\frac{1}{2}%
|\mathbf{y}-\mathbf{y}_{d}|^{2}\,d\mathbf{x}dt+\mathbb{E}\int_{\Gamma
_{T}}\left( \frac{\lambda _{1}}{2}|a|^{2}+\frac{\lambda _{2}}{2}%
|b|^{2}\right) \,d\mathbf{\gamma }dt,  \label{cost}
\end{equation}%
where $\mathbf{y}_{d}\in L_{2}(\Omega \times {\mathcal{O}}_{T})$ is a
desired target field and $\lambda _{1},\lambda _{2}>0.$ We aim to control
the random velocity field $\mathbf{y}$, defined as the solution of the
Stochastic Navier-Stokes equations, through minimization of the cost
functional (\ref{cost}). More precisely, our goal is to solve the following
problem 
\begin{equation*}
(\mathcal{P})\left\{ 
\begin{array}{l}
\underset{(a,b)}{\mbox{minimize}}\{J(a,b,\mathbf{y}):~(a,b)\in \mathcal{A}%
\quad \\ 
\text{and} \\ 
\mathbf{y}\mbox{  is  the weak solution of the  system }\eqref{NSy}%
\mbox{  for
  }(a,b)\in \mathcal{A}\}.%
\end{array}%
\vspace{3mm}\quad \right.
\end{equation*}%
Let us mention that boundary control of fluid flows is of main importance in
several branches of the industry, for instance, in the aviation industry.
The extensive research has been carried out concerning the implementation of
injection-suction devices to control the motion of the fluids (see \cite{bla}%
, \cite{bra}). On the other hand, rotating flow is critically important
across a wide range of scientific, engineering and product applications,
providing design and modeling capability for diverse products such as jet
engines, pumps, food production and vacuum cleaners, as well as geophysical
flows. The control problem for deterministic Newtonian and non-Newtonian
flows, has been widely studied in the literature (see \cite{CC22}, \cite%
{Chem-Cip}, \cite{yas-fer}, \cite{C96}, \cite{ghs}, \cite{gm}). However, it
is well known that the study of turbulent flows, where small random
disturbances produce strong macroscopic effects, requires a statistical
approach. Recently, special attention has been devoted to stochastic optimal
control problems, where control is exerted by a distributed mechanical force
(see \cite{B99}, \cite{CC18}, \cite{CP19}). To the best of our knowledge,
this is the first paper where the boundary control problem is addressed for
stochastic Navier-Stokes equations under Navier-slip boundary conditions.

The plan of the present paper is as follows. In Section \ref{sec1}, we
present the general setting, by introducing the appropriate functional
spaces and some necessary classical inequalities. Section \ref{sec2} deals
with the well-posedness of the state equations. In Section \ref{sec3}, we
show the existence of an optimal solution to the control problem.

\section{General setting}

\label{sec1}\setcounter{equation}{0}

Let $X$ be a real Banach space endowed with the norm $\left\Vert \cdot
\right\Vert _{X}.$ We denote $L_{p}(0,T;X)$ as the space of $X$-valued
measurable $p-$integrable functions defined on $[0,T]$ for $p\geqslant 1$.

For $p,r\geqslant 1$, let $L_{p}(\Omega ,L_{r}(0,T;X))$ be the space of the
processes $\mathbf{v}=\mathbf{v}(\omega ,t)$ with values in $X$ defined on $%
\ \Omega \times \lbrack 0,T],$ adapted to the filtration $\left\{ \mathcal{F}%
_{t}\right\} _{t\in \lbrack 0,T]}$ ,\ and endowed with the norms 
\begin{equation*}
\left\Vert \mathbf{v}\right\Vert _{L_{p}(\Omega ,L_{r}(0,T;X))}=\left( 
\mathbb{E}\left( \int_{0}^{T}\left\Vert \mathbf{v}\right\Vert
_{X}^{r}\,dt\right) ^{\frac{p}{r}}\right) ^{\frac{1}{p}}\text{ }
\end{equation*}%
and%
\begin{equation*}
\left\Vert \mathbf{v}\right\Vert _{L_{p}(\Omega ,L_{\infty }(0,T;X))}=\left( 
\mathbb{E}\sup_{t\in \lbrack 0,T]}\left\Vert \mathbf{v}\right\Vert _{X}^{p}\
\right) ^{\frac{1}{p}}\quad \text{if }r=\infty ,
\end{equation*}%
where $\mathbb{E}$ is the mathematical expectation with respect to the
probability measure $P.$ As usual, in the notation for processes $\mathbf{v}=%
\mathbf{v}(\omega ,t)$, we generally omit the dependence on $\omega \in
\Omega .$

We define the spaces%
\begin{eqnarray*}
H &=&\{\mathbf{v}\in L_{2}({\mathcal{O}}):\,\mbox{div }\mathbf{v}=0\;\ \text{
in}\;\mathcal{D}^{\prime }({\mathcal{O}}),\;\ \mathbf{v}\cdot \mathbf{n}%
=0\;\ \text{ in}\;H^{-1/2}(\Gamma )\}, \\
V &=&\{\mathbf{v}\in H^{1}({\mathcal{O}}):\,\mbox{div }\mathbf{v}=0\;\;\text{%
a.e. in}\;{\mathcal{O}},\;\ \mathbf{v}\cdot \mathbf{n}=0\;\ \text{ in}%
\;H^{1/2}(\Gamma )\}.
\end{eqnarray*}

We denote $(\cdot ,\cdot )$ as the inner product in $L_{2}(\mathcal{O})$ and 
$\Vert \cdot \Vert _{2}$ as the associated norm. The norms in the spaces $%
L_{p}(\mathcal{O})$ and $H^{p}(\mathcal{O})$ are denoted by $\Vert \cdot
\Vert _{p}$ and $\Vert \cdot \Vert _{H^{p}}$. On the space $V$, we consider
the following inner product 
\begin{equation*}
\left( \mathbf{v},\mathbf{z}\right) _{V}=2\left( D\mathbf{v},D\mathbf{z}%
\right) +\alpha \int_{\Gamma }\mathbf{v}\cdot \mathbf{z}
\end{equation*}%
and the corresponding norm $\Vert \mathbf{v}\Vert _{V}=\sqrt{\left( \mathbf{v%
},\mathbf{v}\right) _{V}}.$ \ 

Throughout the article, we often use the continuous embedding results 
\begin{equation}
H^{1}(0,T)\subset C([0,T]),\qquad H^{1}(\mathcal{O})\subset L_{2}(\Gamma ).
\label{aa}
\end{equation}

Let us introduce the notation 
\begin{equation}
\mathbf{v}_{\mathcal{O}}=\int_{\mathcal{O}}\mathbf{v\ }d\mathbf{x.}
\label{000}
\end{equation}%
We notice that for any vector $\mathbf{v}\in V$\ we have $\mathbf{v}_{%
\mathcal{O} }=0,$\ since 
\begin{equation*}
\int_{\mathcal{O} }v_{j}\mathbf{\ }d\mathbf{x}=\int_{\mathcal{O} }\mbox{div }%
(\mathbf{v}x_{j})\mathbf{\ }d\mathbf{x}=\int_{\Gamma }x_{j}(\mathbf{v}\cdot 
\mathbf{n)\ }d\mathbf{\gamma }=0\qquad \text{for}\mathit{\ \ \ }j=1,2.
\end{equation*}%
Using it and the results that can be found on the p. 62, 69 of \cite{lad},
p. 125 of \cite{nir}, and on the p. 16-20 of \cite{tem}, we formulate the
next lemma.

\begin{lemma}
\label{gag} For any $\mathbf{v}\in H^1(\mathcal{O})$ and any $q\geqslant 2$, the
Gagliano-Nirenberg-Sobolev inequality 
\begin{equation}
||\mathbf{v}-\mathbf{v}_{\mathcal{O} }||_{q}\leqslant C||\mathbf{v}%
||_{2}^{2/q}||\nabla \mathbf{v}||_{2}^{1-2/q},  \label{LI}
\end{equation}%
and the trace interpolation inequality%
\begin{equation}
||\mathbf{v}-\mathbf{v}_{\mathcal{O} }||_{L_{q}(\Gamma )}\leqslant C||%
\mathbf{v}||_{2}^{1/q}||\nabla \mathbf{v}||_{2}^{1-1/q}  \label{TT}
\end{equation}%
are valid. Moreover, any $\mathbf{v}\in V$ satisfies Korn's inequality 
\begin{equation}
\left\Vert \mathbf{v}\right\Vert _{H^{1}}\leqslant C\left\Vert \mathbf{v}%
\right\Vert _{V},  \label{korn}
\end{equation}%
that is the norms $\Vert \cdot \Vert _{H^{1}}$ and $\Vert \cdot \Vert _{V}$
are equivalent.
\end{lemma}

\begin{remark}
We should mention that throughout the article, we will represent by $C$ a
generic constant that can assume different values from line to line. These
constants $C$ will depend mainly of the physical constants $\nu ,$ $\alpha , 
$ the domain $\mathcal{O}$, a given time $T>0$.
\end{remark}

Now, we state a formula that can be derived easily via integration by parts 
\begin{equation}
-\int_{\mathcal{O}}\triangle \mathbf{v}\cdot \mathbf{z}\,d\mathbf{x}=2\int_{%
\mathcal{O}}\,D\mathbf{v}\cdot D\mathbf{z}-\int_{\Gamma }2(\mathbf{n}\cdot D%
\mathbf{v})\cdot \mathbf{z},  \label{integrate}
\end{equation}%
which holds for any $\mathbf{v}\in H^{2}(\mathcal{O})$ and $\mathbf{z}\in V$%
. Let us assume that $\mathbf{v}$ satisfies Navier-slip boundary condition %
\eqref{NSy}, then we have 
\begin{equation}
-\int_{\mathcal{O}}\triangle \mathbf{v}\cdot \mathbf{z}\,d\mathbf{x}=(%
\mathbf{v},\mathbf{z}\,)_{V}-\int_{\Gamma }b(\mathbf{z}\cdot {\bm{\tau })}\,d%
\mathbf{\gamma }.  \label{integrate2}
\end{equation}

In what follows we will frequently use 
\begin{equation}
uv\leqslant \varepsilon u^{2}+\frac{v^{2}}{4\varepsilon },\quad \quad
\forall \varepsilon >0,  \label{ab}
\end{equation}%
that is a particular case of Young's inequality 
\begin{equation}
uv\leqslant \frac{u^{p}}{p}+\frac{v^{q}}{q},\quad \quad \frac{1}{p}+\frac{1}{%
q}=1,\quad \forall p,\,q>1.  \label{yi}
\end{equation}

For a vector 
\begin{equation*}
\mathbf{h}=(\mathbf{h}_{1},\dots ,\mathbf{h}_{m})\in H^{m}=\overbrace{%
H\times ...\times H}^{m-times},
\end{equation*}%
we introduce the norm and the absolute value of the inner product of $%
\mathbf{h}$ with a fixed $\mathbf{v}\in H$ as 
\begin{equation}
\left\Vert \mathbf{h}\right\Vert _{2}=\sum_{k=1}^{m}\left\Vert \mathbf{h}%
_{k}\right\Vert _{2}\qquad \text{and}\qquad |\left( \mathbf{h},\mathbf{v}%
\right) |=\left( \sum_{k=1}^{m}\left( \mathbf{h}_{k},\mathbf{v}\right)
^{2}\right) ^{1/2}.  \label{product}
\end{equation}%
Assume that the stochastic noise is represented by 
\begin{equation*}
\mathbf{G}(t,\mathbf{y})\,d{\mathcal{W}}_{t}=\sum_{k=1}^{m}G^{k}(t,\mathbf{y}%
)\,d{\mathcal{W}}_{t}^{k}
\end{equation*}%
where $\mathbf{G}(t,\mathbf{y})=(G^{1}(t,\mathbf{y}),\dots ,G^{m}(t,\mathbf{y%
}))$ has suitable growth assumptions, as defined in the following, and ${%
\mathcal{W}}_{t}=({\mathcal{W}}_{t}^{1},\dots ,{\mathcal{W}}_{t}^{m})$ is a
standard $\mathbb{R}^{m}$-valued Wiener process defined on a complete
probability space $(\Omega ,\mathcal{F},P)$ endowed with a filtration $%
\left\{ \mathcal{F}_{t}\right\} _{t\in \lbrack 0,T]}$. We assume that $%
\mathcal{F}_{0}$\ contains every $P$-null subset of $\Omega $.

Let $\mathbf{G}(t,\mathbf{y}):[0,T]\times H\rightarrow H^{m}$ be Lipschitz
on $\mathbf{y}$ and satisfy the linear growth 
\begin{align}
\left\Vert \mathbf{G}(t,\mathbf{v})-\mathbf{G}(t,\mathbf{z})\right\Vert
_{2}^{2}& \leqslant K\left\Vert \mathbf{v}-\mathbf{z}\right\Vert _{2}^{2}, 
\notag \\
\left\Vert \mathbf{G}(t,\mathbf{v})\right\Vert _{2}& \leqslant K\left(
1+\left\Vert \mathbf{v}\right\Vert _{2}\right) ,\qquad \forall \mathbf{v},%
\mathbf{z}\in H,\;t\in \lbrack 0,T]  \label{G}
\end{align}%
for some positive constant $K.$

Let us define the space of functions \ $\mathcal{H}_{p}(\Gamma )=\left\{
(a,b):||(a,b)||_{\mathcal{H}_{p}(\Gamma )}<+\infty \right\} $ \ with the norm%
\begin{equation*}
||(a,b)||_{\mathcal{H}_{p}(\Gamma )}=||a||_{W_{p}^{1-\frac{1}{p}}(\Gamma
)}+\,||\partial _{t}a||_{H^{\frac{1}{2}}(\Gamma )}+\Vert b\Vert _{W_{p}^{-%
\frac{1}{p}}(\Gamma )}+\Vert b\Vert _{L_{2}(\Gamma )}+||\partial
_{t}b||_{H^{-\frac{1}{2}}(\Gamma )}.
\end{equation*}%
In this work, we consider the data $a,b$ and $\mathbf{u}_{0}$\ belong to the
following Banach spaces 
\begin{equation}
(a,b)\in L_{2}(\Omega \times (0,T);\mathcal{H}_{p}(\Gamma ))\quad \text{for
given}\;p\in (2,+\infty ),\qquad \mathbf{u}_{0}\in L_{2}(\Omega ;H).
\label{eq00sec12}
\end{equation}
In addition, we assume that $(a,b)$ is a pair of predictable stochastic
processes.

%%%%%%%%%%%%%%%%%%%%%%%%%%%%%%%%%%%%%%%%%%%%%%%%%%%%%%%

%%%%%%%%%%%%%%%%%%%%%%%%%%%%%%%%%%%%%%%%%%%%%%%%%%%%%%%

\section{State equation}

\label{sec2}

\setcounter{equation}{0} 
%%%%%%%%%%%%%%%%%%%%%%%%%%%%%%%%%%%%%%%%%%%%%%%%%%%%%%%
This section is devoted to the study of the state equation. We use the
variatioal approach to show the existence and the uniqueness of solution,
and deduce appropriete estimates to study the control problem.

Since, we are considering non-homogeneous boundary conditions, we first
introduce a suitable change of variables based on the solution of the
non-homogeneous linear Stokes equation, which allows to write the state in
terms of a vector field satisfying a homogeneous Navier-slip boundary
condition.

\begin{lemma}
\label{navier slip} Let $(a,b)$ be a given pair of functions satisfying %
\eqref{eq00sec12}. Then there exists a unique solution 
\begin{equation}
\mathbf{a}\in L_{2}(\Omega ;H^{1}(\left( 0,T\right) \times {\mathcal{O}}%
))\cap L_{2}(\Omega \times (0,T);W_{p}^{1}({\mathcal{O}}))  \label{calderon}
\end{equation}%
of the Stokes problem with the non-homogeneous Navier-slip boundary
condition 
\begin{equation}
\left\{ 
\begin{array}{l}
-\Delta \mathbf{a}+\nabla \pi =0,\quad \nabla \cdot \mathbf{a}=0\quad \text{
in }{\mathcal{O}}, \\ 
\mathbf{a}\cdot \mathbf{n}=a,\quad \left[ 2D(\mathbf{a})\,\mathbf{n}+\alpha 
\mathbf{a}\right] \cdot {\bm{\tau }}=b\quad \text{ on }\Gamma ,%
\end{array}%
\right.  \label{ha}
\end{equation}%
such that%
\begin{equation}
||\mathbf{a}||_{W_{p}^{1}({\mathcal{O}})}+||\partial _{t}\mathbf{a}||_{L_{2}(%
{\mathcal{O}})}\leqslant C||(a,b)||_{\mathcal{H}_{p}(\Gamma )},\quad \text{%
a.e. in }\Omega \times (0,T).  \label{cal}
\end{equation}%
In particular, we have%
\begin{equation*}
\mathbf{a}\in L_{2}(\Omega ;C([0,T];L_{2}({\mathcal{O}})))\cap L_{2}(\Omega
\times (0,T);C(\overline{{\mathcal{O}}})\cap H^{1}({\mathcal{O}})).
\end{equation*}%
{\ }
\end{lemma}

\begin{proof}
Let us introduce the function $\mathbf{c}=\nabla h,$ where $h$ is the unique
solution of the system%
\begin{equation*}
\left\{ 
\begin{array}{l}
-\Delta h=0\quad \text{ in }{\mathcal{O}}, \\ 
\frac{\partial h}{\partial \mathbf{n}}=a\quad \text{ on }\Gamma ,%
\end{array}%
\right. \qquad \text{a.e. in}\quad \Omega \times (0,T),
\end{equation*}%
with $\int_{\Gamma }h\,\,d\mathbf{\gamma }=0.$ Theorem 1.10, p. 15 in \cite%
{gir} implies that the function $\mathbf{c}$ \ satisfies\ the estimates%
\begin{equation}
||\mathbf{c}||_{W_{p}^{1}({\mathcal{O}})}\leqslant C_{p}||a||_{W_{p}^{1-%
\frac{1}{p}}(\Gamma )},\qquad ||\partial _{t}\mathbf{c}||_{H^{1}({\mathcal{O}%
})}\leqslant C||\partial _{t}a||_{W_{2}^{\frac{1}{2}}(\Gamma )},  \label{c}
\end{equation}%
where the constant $C_{p}$ depends on $p$, $2<p<\infty $.

Let us consider the following Stokes problem 
\begin{equation*}
\left\{ 
\begin{array}{l}
-\Delta \mathbf{b}+\nabla \pi =0,\quad \nabla \cdot \mathbf{b}=0\quad \text{
in }{\mathcal{O}}, \\ 
\mathbf{b}\cdot \mathbf{n}=0,\quad \left[ 2D(\mathbf{b})\,\mathbf{n}+\alpha 
\mathbf{b}\right] \cdot {\bm{\tau }}=\widetilde{b}\quad \text{ on }\Gamma ,%
\end{array}%
\right. \qquad \text{a.e. in}\quad \Omega \times (0,T)
\end{equation*}%
with $\widetilde{b}=b-\left[ 2D(\mathbf{\mathbf{c}})\,\mathbf{n}+\alpha 
\mathbf{\mathbf{c}}\right] \cdot {\bm{\tau }\in W_{p}^{-\frac{1}{p}}(\Gamma )%
}$ by \eqref{c} and Lemma 2.4 in \cite{amr2}. Using Theorem 2.1 in \ \cite%
{acev}, we have that there exists a unique solution $\mathbf{b}$ of this
Stokes problem\ such that%
\begin{equation}
||\mathbf{b}||_{W_{p}^{1}({\mathcal{O}})}\leqslant C_{p}||\widetilde{b}%
||_{W_{p}^{-\frac{1}{p}}(\Gamma )},\qquad ||\partial _{t}\mathbf{b}||_{H^{1}(%
{\mathcal{O}})}\leqslant C||\partial _{t}\widetilde{b}||_{W_{2}^{-\frac{1}{2}%
}(\Gamma )}.  \label{b}
\end{equation}%
Due to the regularity \eqref{eq00sec12} and the estimates \eqref{c}-\eqref{b}%
, we conclude that the system \eqref{ha} has the unique solution $\mathbf{%
a=b+c}$, satisfying the first estimate in \eqref{cal}. The second one in %
\eqref{cal} is a direct consequence of the embeedings $W_{2}^{1}(0,T)%
\hookrightarrow C([0,T])$ and $W_{p}^{1}({\mathcal{O}})\hookrightarrow C(%
\overline{{\mathcal{O}}}),$ since $2<p<+\infty .$
\end{proof}

\bigskip

With the help of the solution of the non-homogeneous Stokes equation, we
introduce the notion of solution to the state system (\ref{NSy}).

\begin{definition}
\label{1def} Let the data $(a,b)$ and $\mathbf{u}_{0}$ satisfy the
regularity (\ref{eq00sec12}), and $\mathbf{a}$ be the corresponding solution
of \eqref{ha}. A stochastic process $\mathbf{y}=\mathbf{u}+\mathbf{a}$ with $%
\mathbf{u}\in C([0,T];H)\cap L_{2}(0,T;V),\quad P$-a.e. in $\Omega ,$ \ is a strong (in the
stochastic sense) solution of \eqref{NSy} with $\mathbf{y}_{0}=\mathbf{u}%
_{0}+\mathbf{a}(0)$ \ if \ $P$-a.e. in $\Omega $\ the following equation
holds 
\begin{align}
\left( \mathbf{y}(t),\boldsymbol{\varphi }\right) & =\int_{0}^{t}\left[ -\nu
\left( \mathbf{y},\boldsymbol{\varphi }\right) _{V}+\int_{\Gamma }\nu b(%
\boldsymbol{\varphi }\cdot {\bm{\tau }})\,d\mathbf{\gamma }-((\mathbf{y}%
\cdot \nabla )\mathbf{y},\boldsymbol{\varphi })\,\right] ds  \notag \\
&  \notag \\
\quad & +\left( \mathbf{y}_{0},\boldsymbol{\varphi }\right)
+\int_{0}^{t}\left( \mathbf{G}(s,\mathbf{y}(s)),\boldsymbol{\varphi }\right)
\,d{\mathcal{W}}_{s},\qquad \forall t\in \lbrack 0,T],\quad \forall 
\boldsymbol{\varphi }\in V,  \label{res1}
\end{align}%
where the stochastic integral is defined by 
\begin{equation*}
\int_{0}^{t}\left( \mathbf{G}(s,\mathbf{y}(s)),\boldsymbol{\varphi }\right)
\,d{\mathcal{W}}_{s}=\sum_{k=1}^{m}\int_{0}^{t}\left( G^{k}(s,\mathbf{y}(s)),%
\boldsymbol{\varphi }\right) \,d{\mathcal{W}}_{s}^{k}.
\end{equation*}%
\vspace{2mm}\newline
\end{definition}

\bigskip

The existence of solution for the system \eqref{NSy}-\eqref{ICNS} will be
shown by Galerkin's method. Since the injection operator $I:V\rightarrow H$
is a compact operator, there exists a basis $\{\mathbf{e}_{i}\}\subset V$ of
eigenfunctions verifying the property 
\begin{equation}
\left( \mathbf{v},\mathbf{e}_{i}\right) _{V}=\lambda _{i}\left( \mathbf{v},%
\mathbf{e}_{i}\right) ,\qquad \forall \mathbf{v}\in V,\;i\in \mathbb{N},
\label{y4}
\end{equation}%
which is an orthonormal basis for $H$, and the corresponding sequence $%
\{\lambda _{i}\}$ of eigenvalues verifies $\lambda _{i}>0$, $\forall i\in 
\mathbb{N}$ and $\lambda _{i}\rightarrow \infty $ as $i\rightarrow \infty .$
For the details we refer to Theorem 1, p. 355, of \cite{evans}. Moreover the
ellipticity of the equation (\ref{y4}) and the regularity $\Gamma \in C^{2} $
imply that $\{\mathbf{e}_{i}\}\subset C^{2}({\mathcal{O}})\cap V$.

For any fixed $n\in \mathbb{N}$, we consider the subspace $V_{n}=\mathrm{span%
}\,\{\mathbf{e}_{1},\ldots ,\mathbf{e}_{n}\}$ of $V.$ Taking into account
the relation (\ref{integrate2}), the approximate finite dimensional problem
is: for$\ P$-a.e. in $\Omega \,$ to find $\mathbf{y}_{n}$ in the form 
\begin{equation*}
\mathbf{y}_{n}=\mathbf{u}_{n}+\mathbf{a}\quad \text{with}\quad \mathbf{u}%
_{n}(t)=\sum_{j=1}^{n}\beta _{j}^{n}(t)\mathbf{e}_{j}\quad \text{with }t\in
\lbrack 0,T],
\end{equation*}%
as the solution of the following finite dimensional stochastic differential
equation%
\begin{equation}
\left\{ 
\begin{array}{l}
d\left( \mathbf{y}_{n}\,,\boldsymbol{\varphi }\right) =\left[ -\nu \left( 
\mathbf{y}_{n},\boldsymbol{\varphi }\right) _{V}\,+\nu \int_{\Gamma }b(%
\boldsymbol{\varphi }\cdot {\bm{\tau }})\,d\mathbf{\gamma }-(\left( \mathbf{y%
}_{n}\cdot \nabla )\mathbf{y}_{n},\boldsymbol{\varphi }\right) \,\right] dt%
\vspace{2mm} \\ 
\qquad \qquad \,+\left( \mathbf{G}(t,\mathbf{y}_{n}),\boldsymbol{\varphi }%
\right) \,d{\mathcal{W}}_{t},\qquad \forall t\in (0,T),\quad \forall 
\boldsymbol{\varphi }\in V_{n},\vspace{2mm} \\ 
\mathbf{u}_{n}(0)=\mathbf{u}_{n,0},%
\end{array}%
\right.  \label{y1}
\end{equation}%
where $\mathbf{u}_{n,0}=\sum_{j=1}^{n}\left( \mathbf{u}_{0},\mathbf{e}%
_{j}\right) \mathbf{e}_{j}$ is the orthogonal projection of $\mathbf{u}%
_{0}\in H$ \ into the space $V_{n}.$ From the Parseval's identity we infer
that 
\begin{equation}
\left\Vert \mathbf{u}_{n,0}\right\Vert _{2}\leqslant \left\Vert \mathbf{u}%
_{0}\right\Vert _{2}\quad \text{and}\quad \mathbf{u}_{n,0}\longrightarrow 
\mathbf{u}_{0}\qquad \mbox{ strongly in
	}\ H.\quad  \label{CI}
\end{equation}%
The equation (\ref{y1}) defines a system of $n$ stochastic ordinary
differential equations with locally Lipschitz nonlinearities. Hence, there
exists a local-in-time adapted solution $\mathbf{u}_{n}\in
C([0,T_{n}];V_{n}) $ by classical results \cite{KA}. The next lemma will
establish uniform estimates, which guarantee that $\mathbf{u}_{n}$ is a
global-in-time solution.

\begin{lemma}
\label{existence_state} Let the data $(a,b)$ and $\mathbf{u}_{0}$ satisfy
the regularity (\ref{eq00sec12}). Then the system \eqref{y1} has a solution $%
\mathbf{y}_{n}=\mathbf{u}_{n}+\mathbf{a}$, such that 
\begin{equation*}
\mathbf{u}_{n}\in C([0,T];H)\cap L_{2}(0,T;V),\qquad P\text{-a.e. in 
}\Omega \text{.}
\end{equation*}%
Moreover, there exists a positive constant $C_{0}$, such that for the
function%
\begin{equation}
\xi _{0}(t)=e^{-C_{0}t-C_{0}\int_{0}^{t}||(a,b)||_{\mathcal{H}_{p}(\Gamma
)}^{2}ds},\qquad P\text{-a.e. in }\Omega \text{,}  \label{ksi}
\end{equation}%
and any $t\in \lbrack 0,T]$, the following estimate holds 
\begin{align}
\mathbb{E}\sup_{s\in \lbrack 0,t]}\xi _{0}^{2}(s)\left\Vert \mathbf{u}%
_{n}(s)\right\Vert _{2}^{2}& +\nu \mathbb{E}\int_{0}^{t}\xi
_{0}^{2}(s)\left\Vert \mathbf{u}_{n}\right\Vert _{V}^{2}\,ds  \notag \\
& \leqslant C\left( \mathbb{E}\left\Vert \mathbf{u}_{0}\right\Vert _{2}^{2}+%
\mathbb{E}\int_{0}^{t}\xi _{0}^{2}(s)\text{\textrm{A}}(s)\,ds\right)
\label{ineq1}
\end{align}%
where 
\begin{equation}
\text{\textrm{A}}=||(a,b)||_{\mathcal{H}_{p}(\Gamma )}^{2}+1\in L_{1}(\Omega
\times (0,T))  \label{a}
\end{equation}%
and the positive constants $C_{0}$ and $C$ are independent of the parameter $%
n$, which may depend on the regularity of the boundary $\Gamma $ and  the
physical constants $\nu $ and $\alpha $.
\end{lemma}

\begin{proof}
{\ Let }$\xi _{0}${\ be the function defined by (\ref{ksi}) with a constant }%
$C_{0}${\ to be concretized later on (see expression (\ref{C1}) below). \
For each $n\in \mathbb{N}$, let us set }%
\begin{equation*}
g(t)=\xi _{0}^{2}(t)\left\Vert \mathbf{u}_{n}(t)\right\Vert _{2}^{2}+2\nu
\int_{0}^{t}\xi _{0}^{2}(s)\left\Vert \mathbf{u}_{n}(s)\right\Vert
_{V}^{2}ds,\qquad t\in \lbrack 0,T],
\end{equation*}%
{\ and consider the sequence $\{\tau _{N}^{n}\}$}$_{{N\in \mathbb{N}}}${\ of
the stopping times defined by 
\begin{equation}
\tau _{N}^{n}=\inf \{t\geqslant 0:\quad g(t)\geqslant N\}\wedge T_{n}.
\label{g}
\end{equation}%
}Taking $\boldsymbol{\varphi }=\mathbf{e}_{i}$ for each $i=1,\dots ,n$\ in
the equation (\ref{y1}) and using $\mathbf{y}_{n}=\mathbf{u}_{n}+\mathbf{a}$%
, we obtain 
\begin{align}
d\left( \mathbf{u}_{n}\,,\mathbf{e}_{i}\right) & =[-\nu \left( \mathbf{u}%
_{n}+\mathbf{a},\mathbf{e}_{i}\right) _{V}+\nu \int_{\Gamma }b(\mathbf{e}%
_{i}\cdot {\bm{\tau }})\,d\mathbf{\gamma }\,  \notag \\
& +\left( -\mathbf{\partial }_{t}\mathbf{a}-\left( \left( \mathbf{u}_{n}+%
\mathbf{a}\right) \mathbf{\cdot }\nabla \right) \left( \mathbf{u}_{n}+%
\mathbf{a}\right) ,\mathbf{e}_{i}\right) ]dt\,+\left( \mathbf{G}(t,\mathbf{y}%
_{n}),\mathbf{e}_{i}\right) \,d{\mathcal{W}}_{t}.  \notag
\end{align}

\textit{Step 1. Estimate in the space $H$ up to $\tau _{N}^{n}$}. The It\^{o}
formula gives 
\begin{align*}
d\left( \left( \mathbf{u}_{n}\,,\mathbf{e}_{i}\right) ^{2}\right) & =2\left( 
\mathbf{u}_{n},\mathbf{e}_{i}\right) [-\nu \left( \mathbf{u}_{n},\mathbf{e}%
_{i}\right) _{V}-\nu \left( \mathbf{a},\mathbf{e}_{i}\right) _{V}+\nu
\int_{\Gamma }b(\mathbf{e}_{i}\cdot {\bm{\tau }})\,d\mathbf{\gamma }\, \\
& +\left( -\mathbf{\partial }_{t}\mathbf{a}-\left( \left( \mathbf{u}_{n}+%
\mathbf{a}\right) \mathbf{\cdot }\nabla \right) \left( \mathbf{u}_{n}+%
\mathbf{a}\right) ,\mathbf{e}_{i}\right) ]dt\, \\
& +2\left( \mathbf{u}_{n},\mathbf{e}_{i}\right) \left( \mathbf{G}(t,\mathbf{y%
}_{n}),\mathbf{e}_{i}\right) \,d{\mathcal{W}}_{t}+|\left( \mathbf{G}\left( t,%
\mathbf{y}_{n}\right) ,\mathbf{e}_{i}\right) |^{2}\,dt,
\end{align*}%
where the absolute value in the last term is defined by \eqref{product}.
Summing these equalities over $i=1,\dots ,n,$ we obtain%
\begin{align}
d\left( \left\Vert \mathbf{u}_{n}\right\Vert _{2}^{2}\right) +2\nu
\left\Vert \mathbf{u}_{n}\right\Vert _{V}^{2}dt& =\left[ -2\nu \left( 
\mathbf{a},\mathbf{u}_{n}\right) _{V}+\int_{\Gamma }\left\{ -a(\mathbf{u}%
_{n}\cdot \bm{\tau })^{2}+2\nu b(\mathbf{u}_{n}\cdot \bm{\tau })\right\} \,d%
\mathbf{\ \gamma }\,\right] dt  \notag \\
& -2\left( \partial _{t}\mathbf{a}+\left( \left( \mathbf{u}_{n}+\mathbf{a}%
\right) \mathbf{\cdot }\nabla \right) \mathbf{a},\mathbf{u}_{n}\right) dt\, 
\notag \\
& +\sum_{i=1}^{n}|\left( \mathbf{G}\left( t,\mathbf{y}_{n}\right) ,\mathbf{e}%
_{i}\right) |^{2}\,dt+2\left( \mathbf{G}(t,\mathbf{y}_{n}),\mathbf{u}%
_{n}\right) \,d{\mathcal{W}}_{t}  \notag \\
& =I_{1}dt+I_{2}dt+I_{3}dt+2\left( \mathbf{G}(t,\mathbf{y}_{n}),\mathbf{u}%
_{n}\right) \,d{\mathcal{W}}_{t}.  \label{es1}
\end{align}%
Considering Young's inequality \eqref{ab} for an appropriate $\varepsilon >0$%
, the inequalities \ \eqref{LI}-\eqref{korn} and the regularities %
\eqref{calderon}, \eqref{cal}, we estimate the terms $I_{1},$ $I_{2}$ and $%
I_{3}$. Namely 
\begin{align*}
I_{1}& \leqslant 2\nu ||\mathbf{a}||_{V}||\mathbf{u}_{n}||_{V}+\Vert a\Vert
_{L_{\infty }(\Gamma )}\Vert \mathbf{u}_{n}\Vert _{L_{2}(\Gamma )}^{2}+2\nu
\Vert b\Vert {_{L_{2}(\Gamma )}}\Vert \mathbf{u}_{n}\Vert {_{L_{2}(\Gamma )}}%
\, \\
& \leqslant \frac{C}{\nu }||a||_{W_{p}^{1-\frac{1}{p}}(\Gamma )}^{2}||%
\mathbf{u}_{n}||_{2}^{2}+C\nu (||\mathbf{a}||_{V}^{2}+\Vert b\Vert
_{L_{2}(\Gamma )}^{2})+\frac{\nu }{2}||\mathbf{u}_{n}||_{V}^{2} \\
& \leqslant \frac{C}{\nu }\text{\textrm{A}}||\mathbf{u}_{n}||_{2}^{2}+C\nu
||(a,b)||_{\mathcal{H}_{p}(\Gamma )}^{2}+\frac{\nu }{2}||\mathbf{u}%
_{n}||_{V}^{2}\,,
\end{align*}%
where \textrm{A} is defined by (\ref{a}). A similar reasoning gives%
\begin{eqnarray*}
I_{2} &\leqslant &2\left( \Vert \partial _{t}\mathbf{a}\Vert _{2}+\,||%
\mathbf{a}||_{C(\overline{\Omega })}\Vert \nabla \mathbf{a}\Vert _{2}\right)
\,||\mathbf{u}_{n}||_{2}+2\Vert \nabla \mathbf{a}\Vert _{2}\,\Vert \mathbf{u}%
_{n}\Vert _{4}^{2} \\
&\leqslant &C\left( \Vert \partial _{t}\mathbf{a}\Vert _{2}^{2}+\,||\mathbf{a%
}||_{C(\overline{\Omega })}^{2}+\Vert \nabla \mathbf{a}\Vert
_{2}^{2}+1\right) \left( 1+||\mathbf{u}_{n}||_{2}\right) +\frac{\nu }{2}||%
\mathbf{u}_{n}||_{V}^{2} \\
&\leqslant &C\text{\textrm{A}}(1+\,||\mathbf{u}_{n}||_{2}^{2})\,+\frac{\nu }{%
2}||\mathbf{u}_{n}||_{V}^{2}
\end{eqnarray*}%
and 
\begin{eqnarray*}
I_{3} &=&\sum_{i=1}^{n}|\left( \mathbf{G}\left( t,\mathbf{y}_{n}\right) ,%
\mathbf{e}_{i}\right) |^{2}\leqslant C||\mathbf{G}\left( t,\mathbf{y}%
_{n}\right) ||_{2}^{2}\leqslant C(1+||\mathbf{y}_{n}||_{2}^{2}) \\
&\leqslant &C(1+||\mathbf{u}_{n}||_{2}^{2}+||\mathbf{a}||_{2}^{2})\leqslant
C(||\mathbf{u}_{n}||_{2}^{2}+\text{\textrm{A}}),
\end{eqnarray*}%
where we used the assumption (\ref{G}). Gathering the previous estimates, we
obtain the existence of a positive constant $C_{0},$ such that 
\begin{equation}
I_{1}+I_{2}+I_{3}\leqslant 2C_{0}\text{\textrm{A}}(||\mathbf{u}_{n}\Vert
_{2}^{2}+1).  \label{C1}
\end{equation}

Taking the function $\xi _{0}$ as in \eqref{ksi}, thanks to \eqref{es1}-%
\eqref{C1}, the application of It\^{o}'s formula yields 
\begin{align*}
\xi _{0}^{2}(s)\left\Vert \mathbf{u}_{n}(s)\right\Vert _{2}^{2}& +2\nu
\int_{0}^{s}\xi _{0}^{2}(r)\left\Vert \mathbf{u}_{n}\right\Vert _{V}^{2}dr \\
& =\left\Vert \mathbf{u}_{n}(0)\right\Vert _{2}^{2}-2C_{0}\int_{0}^{s}\xi
_{0}^{2}(r)\text{\textrm{A}}(r)\left\Vert \mathbf{u}_{n}\right\Vert
_{2}^{2}dr \\
& +\int_{0}^{s}\xi _{0}^{2}(r)(I_{1}+I_{2}+I_{3})dr+2\int_{0}^{s}\xi
_{0}^{2}(r)\left( \mathbf{G}(r,\mathbf{y}_{n}),\mathbf{u}_{n}\right) \,d{%
\mathcal{W}}_{r} \\
& \leqslant \left\Vert \mathbf{u}_{n,0}\right\Vert
_{2}^{2}-2C_{0}\int_{0}^{s}\xi _{0}^{2}(r)\text{\textrm{A}}(r)\left\Vert 
\mathbf{u}_{n}\right\Vert _{2}^{2}dr \\
& +2C_{0}\int_{0}^{s}\xi _{0}^{2}(r)\text{\textrm{A}}(r)dr+2\int_{0}^{s}\xi
_{0}^{2}(r)\left( \mathbf{G}(r,\mathbf{y}_{n}),\mathbf{u}_{n}\right) \,d{%
\mathcal{W}}_{r} \\
& \leqslant \left\Vert \mathbf{u}_{n,0}\right\Vert
_{2}^{2}+2C_{0}\int_{0}^{s}\xi _{0}^{2}(r)\text{\textrm{A}}%
(r)dr+2\int_{0}^{s}\xi _{0}^{2}(r)\left( \mathbf{G}(r,\mathbf{y}_{n}),%
\mathbf{u}_{n}\right) \,d{\mathcal{W}}_{r}.
\end{align*}%
Therefore, we can write 
\begin{align}
\xi _{0}^{2}(s)\left\Vert \mathbf{u}_{n}(s)\right\Vert _{2}^{2}+\nu
\int_{0}^{s}\xi _{0}^{2}(r)\left\Vert \mathbf{u}_{n}\right\Vert _{V}^{2}\
dr& \leqslant \left\Vert \mathbf{u}_{n,0}\right\Vert
_{2}^{2}+C\int_{0}^{s}\xi _{0}^{2}(r)A(r)\,dr  \notag \\
& +2\int_{0}^{s}\xi _{0}^{2}(r)\left( \mathbf{G}\left( r,\mathbf{y}%
_{n}\right) ,\mathbf{u}_{n}\right) \,d{\mathcal{W}}_{r}.  \label{y2}
\end{align}%
Now, considering the sequence $(\tau _{N}^{n})$ of the stopping times
introduced in \eqref{g} and using \eqref{G}, the Burkholder-Davis-Gundy
inequality gives 
\begin{align*}
& \mathbb{E}\sup_{s\in \lbrack 0,\tau _{N}^{n}\wedge t]}\left\vert
\int_{0}^{s}\xi _{0}^{2}(r)\left( \mathbf{G}\left( r,\mathbf{y}_{n}\right) ,%
\mathbf{u}_{n}\right) \,d{\mathcal{W}}_{r}\right\vert \leqslant \mathbb{E}%
\left( \int_{0}^{\tau _{N}^{n}\wedge t}\xi _{0}^{4}(s)\left\vert \left( 
\mathbf{G}\left( s,\mathbf{y}_{n}\right) ,\mathbf{u}_{n}\right) \right\vert
^{2}\,ds\right) ^{\frac{1}{2}} \\
& \qquad \qquad \leqslant \mathbb{E}\sup_{s\in \lbrack 0,\tau _{N}^{n}\wedge
t]}\xi _{0}(s)\left\Vert \mathbf{u}_{n}(s)\right\Vert _{2}\left(
\int_{0}^{\tau _{N}^{n}\wedge t}\xi _{0}^{2}(s)\left\Vert \mathbf{G}\left( s,%
\mathbf{y}_{n}\right) \right\Vert _{2}^{2}\,ds\right) ^{\frac{1}{2}} \\
& \qquad \qquad \leqslant \varepsilon \,\mathbb{E}\sup_{s\in \lbrack 0,\tau
_{N}^{n}\wedge t]}\xi _{0}^{2}(s)\Vert \mathbf{u}_{n}(s)\Vert
_{2}^{2}+C_{\varepsilon }\mathbb{E}\int_{0}^{\tau _{N}^{n}\wedge t}\xi
_{0}^{2}(s)(||\mathbf{u}_{n}||_{2}^{2}+\text{\textrm{A}}(s))\,ds.
\end{align*}%
For $t\in \lbrack 0,T]$, we first take the supremun of the relation %
\eqref{y2} for $s\in \lbrack 0,\tau _{N}^{n}\wedge t]$, next we take the
expectation and incorporate the previous estimate of the stochastic term
with $\varepsilon =\frac{1}{2}$. Then considering (\ref{CI}), we deduce 
\begin{align*}
& \frac{1}{2}\mathbb{E}\sup_{s\in \lbrack 0,\tau _{N}^{n}\wedge t]}\xi
_{0}^{2}(s)\Vert \mathbf{u}_{n}(s)\Vert _{2}^{2}+\nu \mathbb{E}%
\int_{0}^{\tau _{N}^{n}\wedge t}\xi _{0}^{2}(s)\left\Vert \mathbf{u}%
_{n}\right\Vert _{V}^{2}\,ds \\
& \qquad \qquad \leqslant \mathbb{E}\left\Vert \mathbf{u}_{0}\right\Vert
_{2}^{2}+C\mathbb{E}\int_{0}^{\tau _{N}^{n}\wedge t}\xi _{0}^{2}(s)\text{%
\textrm{A}}(s)\,ds+C\mathbb{E}\int_{0}^{\tau _{N}^{n}\wedge t}\xi
_{0}^{2}\left\Vert \mathbf{u}_{n}\right\Vert _{2}^{2}\,ds.
\end{align*}%
Hence, the function 
\begin{equation*}
f(t)=\mathbb{E}\sup_{s\in \lbrack 0,\tau _{N}^{n}\wedge t]}\xi
_{0}^{2}(s)\Vert \mathbf{u}_{n}(s)\Vert _{2}^{2}+2\nu \mathbb{E}%
\int_{0}^{\tau _{N}^{n}\wedge t}\xi _{0}^{2}(s)\left\Vert \mathbf{u}%
_{n}\right\Vert _{V}^{2}\,ds
\end{equation*}%
fulfills the Gronwall type inequality 
\begin{equation*}
\frac{1}{2}f(t)\leqslant \mathbb{E}\left\Vert \mathbf{u}_{n}(0)\right\Vert
_{2}^{2}+C\mathbb{E}\int_{0}^{t}\xi _{0}^{2}(s)\text{\textrm{A}}%
(s)\,ds+\int_{0}^{t}f(s)ds,
\end{equation*}%
which implies 
\begin{align}
\mathbb{E}\sup_{s\in \lbrack 0,\tau _{N}^{n}\wedge t]}\xi
_{0}^{2}(s)\left\Vert \mathbf{u}_{n}(s)\right\Vert _{2}^{2}& +2\nu \mathbb{E}%
\int_{0}^{\tau _{N}^{n}\wedge t}\xi _{0}^{2}(s)\left\Vert \mathbf{u}%
_{n}\right\Vert _{V}^{2}\,ds  \notag \\
& \leqslant C\mathbb{E}\left\Vert \mathbf{u}_{0}\right\Vert _{2}^{2}+C%
\mathbb{E}\int_{0}^{t}\xi _{0}^{2}(s)\text{\textrm{A}}(s)\,ds.  \label{IMP_1}
\end{align}%
\textit{Step 2. The limit transition as $N\rightarrow \infty $.} From (\ref%
{IMP_1}) we have 
\begin{equation*}
\mathbb{E}\sup_{s\in \lbrack 0,\tau _{N}^{n}\wedge T]}g(s)\leqslant C
\end{equation*}%
for some constant $C$ independent of $N$ and $n$. Let us fix $n\in \mathbb{N}
$. Since $\mathbf{u}_{n}\in C([0,T_{n}];V_{n}),$ we have $\ g(\tau
_{N}^{n})\geqslant N$ and%
\begin{align}
\mathbb{E}\sup_{s\in \lbrack 0,\tau _{N}^{n}\wedge T]}g(s)& \geqslant 
\mathbb{E}\left( \sup_{s\in \lbrack 0,\tau _{N}^{n}\wedge T]}1_{\{\tau
_{N}^{n}<T\}}g(s)\right)  \notag \\
& =\mathbb{E}\left( 1_{\{\tau _{N}^{n}<T\}}\ g(\tau _{N}^{n})\right)
\geqslant NP\left( \tau _{N}^{n}<T\right) ,  \label{fc}
\end{align}%
which implies that $P\left( \tau _{N}^{n}<T\right) \rightarrow 0$, as $%
N\rightarrow \infty .$ This means that $\tau _{N}^{n}\rightarrow T$ in
probability as $N\rightarrow \infty $. Then, there exists a subsequence $%
\{\tau _{N_{k}}^{n}\}$ of $\{\tau _{N}^{n}\}$ (which may depend on $n$) such
that 
\begin{equation*}
\tau _{N_{k}}^{n}(\omega )\rightarrow T\text{\qquad for a. e. \ }\omega \in
\Omega \quad \text{as }k\rightarrow \infty .
\end{equation*}%
Since $\tau _{N_{k}}^{n}\leqslant T_{n}\leqslant T$, we deduce that $T_{n}=T$%
, hence $\mathbf{y}_{n}=\mathbf{u}_{n}+\mathbf{a}$ is a global-in-time
solution of the stochastic differential equation (\ref{y1}). In addition,
for each fixed $n\in \mathbb{N}$, the sequence $\left\{ \tau
_{N}^{n}\right\} $ is monotone on $N$, therefore we can apply the monotone
convergence theorem in order to pass to the limit in the inequality (\ref%
{IMP_1}) as $N\rightarrow \infty $, \ thereby deducing the estimate (\ref%
{ineq1}).
\end{proof}

\bigskip

In the next lemma, by assuming a better integrability for the initial data,
we improve the integrability properties for the solution $\mathbf{y}_{n}$ of
problem (\ref{y1}).

\begin{lemma}
\label{existence_state1} Let the data $(a,b)$ and $\mathbf{u}_{0}$ satisfy
the regularity (\ref{eq00sec12}). In addition we assume 
\begin{eqnarray}
(a,b) &\in &L_{4}(\Omega \times (0,T);\mathcal{H}_{p}(\Gamma )),  \notag \\
\mathbf{u}_{0} &\in &L_{4}(\Omega ;H).  \label{0}
\end{eqnarray}%
Then, the solution $\mathbf{y}_{n}=\mathbf{u}_{n}+\mathbf{a}$ of problem $(%
\ref{y1})$ has the regularity 
\begin{equation*}
\mathbf{u}_{n}\in C([0,T];H)\cap L_{4}(0,T;V),\quad P\text{-a.e. in }%
\Omega \text{,}
\end{equation*}%
such that 
\begin{eqnarray}
&&\mathbb{E}\sup_{s\in \lbrack 0,t]}\xi _{0}^{4}(s)\left\Vert \mathbf{u}%
_{n}(s)\right\Vert _{2}^{4}+8\nu ^{2}\mathbb{E}\left( \int_{0}^{t}\xi
_{0}^{2}(s)\left\Vert \mathbf{u}_{n}(s)\right\Vert _{V}^{2}ds\right) ^{2} 
\notag \\
&&\qquad \qquad \qquad \qquad \qquad \leqslant C\left( \mathbb{E}\left\Vert 
\mathbf{u}_{0}\right\Vert _{2}^{4}+\mathbb{E}\int_{0}^{t}\xi _{0}^{4}\text{%
\textrm{B}}(s)\,ds\right) ,\quad t\in \lbrack 0,T],  \label{lp1}
\end{eqnarray}%
where the function $\xi _{0}$ is defined in (\ref{ksi}), 
\begin{equation}
\text{\textrm{B}}=||(a,b)||_{\mathcal{H}_{p}(\Gamma )}^{4}+1\in L_{1}(\Omega
\times (0,T)),  \label{BB}
\end{equation}%
and $C$ is a positive constant, being independent of $n$.
\end{lemma}

\begin{proof}
Taking the square on both sides of the inequality (\ref{y2}) and the
supremum on $s\in \lbrack 0,\tau _{N}^{n}\wedge t]$ \ with $\tau _{N}^{n}$
defined by (\ref{g}), we infer that 
\begin{align*}
\sup_{s\in \lbrack 0,\tau _{N}^{n}\wedge t]}\xi _{0}^{4}(s)\left\Vert 
\mathbf{u}_{n}(s)\right\Vert _{2}^{4}& +\nu ^{2}\left( \int_{0}^{\tau
_{N}^{n}\wedge t}\xi _{0}^{2}(s)\left\Vert \mathbf{u}_{n}(s)\right\Vert
_{V}^{2}ds\right) ^{2} \\
& \leqslant 8\left( \left\Vert \mathbf{u}_{n,0}\right\Vert _{2}^{4}+C^{2}%
\mathbb{E}\int_{0}^{\tau _{N}^{n}\wedge t}\xi _{0}^{4}(s)\text{\textrm{B}}%
(s)\,ds\right) \\
& +4\sup_{s\in \lbrack 0,\tau _{N}^{n}\wedge t]}\left\vert \int_{0}^{s}\xi
_{0}^{2}(r)\left( \mathbf{G}\left( r,\mathbf{y}_{n}\right) ,\mathbf{u}%
_{n}\right) \,d{\mathcal{W}}_{r}\right\vert ^{2}
\end{align*}%
where \textrm{B} is defined by (\ref{BB}). Therefore taking the expectation
in this inequality and applying the Burkholder-Davis-Gundy inequality 
\begin{align*}
\mathbb{E}& \sup_{s\in \lbrack 0,\tau _{N}^{n}\wedge t]}|\int_{0}^{s}\xi
_{0}^{2}(r)\left( \mathbf{G}\left( r,\mathbf{y}_{n}\right) ,\mathbf{u}%
_{n}\right) \,d{\mathcal{W}}_{r}|^{2}\leqslant \mathbb{E}\left(
\int_{0}^{\tau _{N}^{n}\wedge t}\xi _{0}^{4}(s)\left\vert \left( \mathbf{G}%
\left( s,\mathbf{\ y}_{n}\right) ,\mathbf{u}_{n}\right) \right\vert
^{2}\,ds\right) \\
& \qquad \qquad \leqslant \mathbb{E}\sup_{s\in \lbrack 0,\tau _{N}^{n}\wedge
t]}\xi _{0}^{2}\left\Vert \mathbf{u}_{n}\right\Vert _{2}^{2}\int_{0}^{\tau
_{N}^{n}\wedge t}\xi _{0}^{2}\left\Vert \mathbf{G}\left( s,\mathbf{y}%
_{n}\right) \right\Vert _{2}^{2}\,ds \\
& \qquad \qquad \leqslant \frac{1}{2}\,\mathbb{E}\sup_{s\in \lbrack 0,\tau
_{N}^{n}\wedge t]}\xi _{0}^{4}(s)\Vert \mathbf{u}_{n}(s)\Vert _{2}^{4}+C%
\mathbb{E}\int_{0}^{\tau _{N}^{n}\wedge t}\xi _{0}^{4}(s)\left( ||\mathbf{u}%
_{n}||_{2}^{4}+\text{\textrm{B}}(s)\right) \,ds,
\end{align*}%
we obtain 
\begin{align*}
\frac{1}{2}\mathbb{E}& \sup_{s\in \lbrack 0,\tau _{N}^{n}\wedge t]}\xi
_{0}^{4}(s)\left\Vert \mathbf{u}_{n}(s)\right\Vert _{2}^{4}+\nu ^{2}\mathbb{E%
}\left( \int_{0}^{\tau _{N}^{n}\wedge t}\xi _{0}^{2}(s)\left\Vert \mathbf{u}%
_{n}\right\Vert _{V}^{2}ds\right) ^{2} \\
& \qquad \leqslant C(\mathbb{E}\left\Vert \mathbf{u}_{0}\right\Vert
_{2}^{4}+\int_{0}^{\tau _{N}^{n}\wedge t}\xi _{0}^{4}(s)\text{\textrm{B}}%
(s)\,ds+C\mathbb{E}\int_{0}^{\tau _{N}^{n}\wedge t}\xi _{0}^{4}(s)(1+||%
\mathbf{u}_{n}||_{2}^{4})\,ds.
\end{align*}%
Using Gronwall's inequality, we deduce that 
\begin{eqnarray}
\mathbb{E}\sup_{s\in \lbrack 0,\tau _{N}^{n}\wedge t]}\xi
_{0}^{4}(s)\left\Vert \mathbf{u}_{n}(s)\right\Vert _{2}^{4} &+&\nu ^{2}%
\mathbb{E}\left( \int_{0}^{\tau _{N}^{n}\wedge t}\xi _{0}^{2}(s)\left\Vert 
\mathbf{u}_{n}\right\Vert _{V}^{2}ds\right) ^{2}\leqslant C\mathbb{E}%
\left\Vert \mathbf{u}_{0}\right\Vert _{2}^{4}  \notag \\
&+&C\int_{0}^{\tau _{N}^{n}\wedge t}\xi _{0}^{4}(s)\text{\textrm{B}}%
(s)\,ds,\quad \forall t\in \lbrack 0,T].  \label{IMP_7}
\end{eqnarray}%
Arguing as in the proof of Lemma \ref{existence_state}, there exists a
monotone subsequence $\{\tau _{N_{k}}^{n}\}$ of $\{\tau _{N}^{n}\}$, which
converges to $T$ a.e. $\omega \in \Omega $, as $k\rightarrow \infty $. Thus,
applying the monotone convergence theorem, we can pass to the limit in (\ref%
{IMP_7}) as $k\rightarrow \infty $, in order to deduce the estimate (\ref%
{lp1}).
\end{proof}

\bigskip

\begin{theorem}
\label{the_1} Let the data $(a,b)$ and $\mathbf{u}_{0}$ satisfy the
regularity (\ref{eq00sec12}) and (\ref{0}). Then there exists, a unique
strong solution $\mathbf{y}=\mathbf{u}+\mathbf{a}$ to the system \eqref{NSy}-%
\eqref{ICNS}, such that 
\begin{equation*}
\mathbf{u}\in C([0,T];H)\cap L_{4}(0,T;V),\quad P\text{-a.e. in }%
\Omega ,
\end{equation*}%
and for any $t\in \lbrack 0,T],$ the following estimates hold 
\begin{align}
\mathbb{E}\sup_{s\in \lbrack 0,t]}\xi _{0}^{2}(s)\left\Vert \mathbf{u}%
(s)\right\Vert _{2}^{2}& +\nu \mathbb{E}\int_{0}^{t}\xi
_{0}^{2}(s)\left\Vert \mathbf{u}\right\Vert _{V}^{2}\,ds  \notag \\
& \leqslant C\left( \mathbb{E}\left\Vert \mathbf{u}_{0}\right\Vert _{2}^{2}+%
\mathbb{E}\int_{0}^{t}\xi _{0}^{2}(s)\mathrm{A}(s)\,ds\right) ,
\label{21:47}
\end{align}%
\begin{align}
\mathbb{E}\sup_{s\in \lbrack 0,t]}\xi _{0}^{4}(s)\left\Vert \mathbf{u}%
(s)\right\Vert _{2}^{4}& +\nu ^{2}\mathbb{E}\left( \int_{0}^{t}\xi
_{0}^{2}(s)\left\Vert \mathbf{u}\right\Vert _{V}^{2}ds\right) ^{2} 
\notag \\
& \leqslant C\left( \mathbb{E}\left\Vert \mathbf{y}_{0}\right\Vert
_{2}^{4}+\nu ^{2}\mathbb{E}\int_{0}^{t}\xi _{0}^{4}(s)\mathrm{B}%
(s)\,ds\right) ,  \label{uny}
\end{align}%
where the functions $\xi _{0}$ and \textrm{A}\textsl{, }\textrm{B} are
defined by (\ref{ksi}) and (\ref{a}), (\ref{BB}), respectively. Here $C$ is
a positive constant that is independent of $n$.\ 
\end{theorem}

\begin{proof}
The proof is splitted into three steps. \vspace{2mm} \newline
\textit{Step 1. Convergence related to the projection operator.} Let $%
P_{n}:V\rightarrow V_{n}$ be the orthogonal projection defined by%
\begin{equation*}
P_{n}\mathbf{v}=\sum_{j=1}^{n}\widetilde{\beta }_{j}\widetilde{\mathbf{e}}%
_{j}=\sum_{j=1}^{n}\beta _{j}\mathbf{e}_{j}\quad \text{with }\widetilde{%
\beta }_{j}=\left( \mathbf{v},\widetilde{\mathbf{e}}_{j}\right) _{V}\quad 
\text{and}\quad \beta _{j}=\left( \mathbf{v},\mathbf{e}_{j}\right) ,\quad
\forall \mathbf{v}\in V,
\end{equation*}%
where $\{\widetilde{\mathbf{e}}_{j}=\frac{1}{\sqrt{\lambda _{j}}}\mathbf{e}%
_{j}\}_{j=1}^{\infty }$ is the orthonormal basis of $V.$ \ By Parseval's
identity, for any $\mathbf{v}\in V$ we have 
\begin{eqnarray}
&||P_{n}\mathbf{v}||_{2}\leqslant &||\mathbf{v}||_{2},\qquad ||P_{n}\mathbf{v%
}||_{V}\leqslant ||\mathbf{v}||_{V},  \notag \\
&&\qquad P_{n}\mathbf{v}\longrightarrow \mathbf{v}\qquad 
\mbox{ strongly in
	}\text{ }V.  \label{2}
\end{eqnarray}%
Considering an arbitrary $\mathbf{z}\in L_{s}(\Omega \times (0,T);V)$ for
some $s\geqslant 1$, we have%
\begin{equation*}
||P_{n}\mathbf{z}||_{V}\leqslant ||\mathbf{z}||_{V}\quad \mbox{ and
	}\quad P_{n}\mathbf{z}(\omega ,t)\rightarrow \mathbf{z}(\omega ,t)\qquad 
\text{ strongly in}\ V,
\end{equation*}%
which are valid $P$-a.e. $\omega \in \Omega $ and a.e. $t\in (0,T).$ Hence,
Lebesgue's dominated convergence theorem implies that for any $\mathbf{z}\in
L_{s}(\Omega \times (0,T);V)$, we have%
\begin{equation}
P_{n}\mathbf{z}\longrightarrow \mathbf{z}\qquad \mbox{ strongly in
	}\ L_{s}(\Omega \times (0,T);V).  \label{c02}
\end{equation}%
\textit{Step 2.} \textit{Passage to the limit in the weak sense.}

Let us define$~f_{0}(t)=C_{0}\left(||(a,b)||_{\mathcal{H}_{p}(\Gamma
)}^{2}+1\right)$. Since 
\begin{equation*}
\int_{0}^{T}f_{0}(s)\ ds\leqslant C(\omega )<+\infty \quad \text{for all }%
\omega \in \Omega \backslash A,\quad \text{where }P(A)=0
\end{equation*}%
by \ {(\ref{eq00sec12}), }there exists a positive constant $K(\omega ),$
which depends only on $\omega \in \Omega \backslash A$ and satisfies%
\begin{equation}
0<K(\omega )\leqslant \xi _{0}(t)=e^{-\int_{0}^{t}f_{0}(s)\ ds}\leqslant
1\quad \text{for all }\omega \in \Omega \backslash A,\quad t\in \lbrack 0,T].
\label{q2}
\end{equation}

The estimates (\ref{ineq1}) and (\ref{lp1}) give that%
\begin{eqnarray}
\mathbb{E}\sup_{t\in \lbrack 0,T]}\left\Vert \xi _{0}(t)\mathbf{u}%
_{n}(t)\right\Vert _{2}^{2} &\leqslant &C,\qquad \mathbb{E}%
\int_{0}^{T}\left\Vert \xi _{0}\mathbf{u}_{n}\right\Vert _{V}^{2}\
dt\leqslant C,  \notag \\
\mathbb{E}\sup_{t\in \lbrack 0,T]}\left\Vert \xi _{0}(t)\mathbf{u}%
_{n}(t)\right\Vert _{2}^{4} &\leqslant &C,\qquad \mathbb{E}\left(
\int_{0}^{T}\left\Vert \xi _{0}\mathbf{u}_{n}\right\Vert _{V}^{2}\ dt\right)
^{2}\leqslant C  \label{as}
\end{eqnarray}%
for some constant $C$ that is independent of the index $n$. These uniform estimates
imply 
\begin{equation}
\left\Vert \xi _{0}^{2}\left( \mathbf{y}_{n}\cdot \nabla \right) \mathbf{y}%
_{n}\right\Vert _{L_{2}(\Omega \times (0,T);V^{\prime })}\leqslant C,\quad
\forall n\in \mathbb{N},  \label{l01}
\end{equation}%
where $V^{\prime }$ denotes the topological dual of the space $V.$ The
uniform estimates (\ref{as}) ensures the existence of a suitable subsequence 
$\mathbf{u}_{n}$, which is indexed by the same index $n$ to simplify the
notation, and a function $\mathbf{u},$ such that 
\begin{eqnarray}
\xi _{0}\mathbf{u}_{n} &\rightharpoonup &\xi _{0}\mathbf{u}\quad 
\mbox{ weakly in
}\ L_{2}(\Omega \times (0,T);V)\cap L_{4}(\Omega ,L_{2}(0,T;V)),  \notag \\
\xi _{0}\mathbf{u}_{n} &\rightharpoonup &\xi _{0}\mathbf{u}\quad 
\mbox{ *-weakly in
}\ L_{2}(\Omega ,L_{\infty }(0,T;H))\cap L_{4}(\Omega ,L_{\infty
}(0,T;H)).\quad  \label{c1}
\end{eqnarray}%
Moreover, we have 
\begin{equation}
\xi _{0}P_{n}\mathbf{u}\longrightarrow \xi _{0}\mathbf{u}\qquad 
\mbox{ strongly in
}\ L_{2}(\Omega \times (0,T);V)\cap L_{4}(\Omega ,L_{2}(0,T;V))  \label{c02Y}
\end{equation}%
by (\ref{c02}). The limit function $\mathbf{u}$ satisfies the estimates %
\eqref{21:47}, \eqref{uny} by the lower semicontinuity of integral in $L_{2}$
and $L_{4}$ spaces.

Considering (\ref{G}) and (\ref{l01}), there exist some operators $B^{\ast
}(t)$ and $\mathbf{G}^{\ast }(t)$ such that%
\begin{eqnarray}
\xi _{0}\mathbf{G}(t,\mathbf{y}_{n}) &\rightharpoonup &\xi _{0}\mathbf{G}%
^{\ast }(t)\qquad \mbox{ weakly in
}\ L_{2}(\Omega \times (0,T);H^{m}),  \notag \\
\xi _{0}^{2}\left( \mathbf{y}_{n}\cdot \nabla \right) \mathbf{y}_{n}
&\rightharpoonup &\xi _{0}^{2}B^{\ast }(t)\qquad \mbox{ weakly in
}\ L_{2}(\Omega \times (0,T); V^{\prime }).  \label{c01}
\end{eqnarray}

Since $\mathbf{y}_{n}$ solves the equation (\ref{y1}), then using It\^{o}'s
formula, we infer that 
\begin{align*}
d\left( \xi _{0}^{2}\mathbf{y}_{n},\boldsymbol{\varphi }\right) & =\xi
_{0}^{2}[-\nu \left( \mathbf{y}_{n},\boldsymbol{\varphi }\right) _{V}\,+\nu
\int_{\Gamma }b(\boldsymbol{\varphi }\cdot {\bm{\tau }})\,d\mathbf{\gamma }%
-(\left( \mathbf{y}_{n}\cdot \nabla )\mathbf{y}_{n},\boldsymbol{\varphi }%
\right) \, \\
& -2f_{0}(t)\left( \mathbf{y}_{n},\boldsymbol{\varphi }\right) ]\,dt+\xi
_{0}^{2}\left( \mathbf{G}(t,\mathbf{y}_{n}),\boldsymbol{\varphi }\right) \,d{%
\mathcal{W}}_{t},
\end{align*}%
that is, the following integral equation holds 
\begin{align}
\left( \xi _{0}^{2}(t)\mathbf{y}_{n}(t),\boldsymbol{\varphi }\right)
&-\left( \mathbf{y}_{n,0},\boldsymbol{\varphi }\right) =\int_{0}^{t}\xi
_{0}^{2}(s)[-\nu \left( \mathbf{y}_{n}(s),\boldsymbol{\varphi }\right)
_{V}\,+\nu \int_{\Gamma }b(s)(\boldsymbol{\varphi }\cdot {\bm{\tau }})\,d%
\mathbf{\gamma }  \notag \\
&-(\left( \mathbf{y}_{n}(s)\cdot \nabla )\mathbf{y}_{n}(s),\boldsymbol{%
\varphi }\right) \, -2f_{0}(s)\left( \mathbf{y}_{n}(s),\boldsymbol{\varphi }%
\right) ]\,ds  \notag \\
&+\int_{0}^{t}\xi _{0}^{2}(s)\left( \mathbf{G}(s,\mathbf{y}_{n}),\boldsymbol{%
\varphi }\right) \,d{\mathcal{W}}_{s},\quad \forall t\in \lbrack 0,T],\quad P%
\text{-a.e.in }\Omega .  \label{CC1}
\end{align}

Denoting 
\begin{equation*}
\mathbf{h}_{n}(t)=\xi _{0}^{2}(t)\mathbf{y}_{n}(t)-\int_{0}^{t}\xi
_{0}^{2}(s)\mathbf{G}(s,\mathbf{y}_{n})\,d{\mathcal{W}}_{s}
\end{equation*}%
the following differential equation holds 
\begin{align}
\frac{\partial }{\partial t}\left( \mathbf{h}_{n}(t),\boldsymbol{\varphi }%
\right) & =\xi _{0}^{2}(t)[-\nu \left( \mathbf{y}_{n}(t),\boldsymbol{\varphi 
}\right) _{V}\,+\nu \int_{\Gamma }b(t)(\boldsymbol{\varphi }\cdot {\bm{\tau }%
})\,d\mathbf{\gamma }-(\left( \mathbf{y}_{n}(t)\cdot \nabla )\mathbf{y}%
_{n}(t),\boldsymbol{\varphi }\right) \,  \notag  \label{11:33} \\
& -2f_{0}(t)\left( \mathbf{y}_{n}(t),\boldsymbol{\varphi }\right) ],\quad P%
\text{-a.e. in }\Omega ,\quad \forall t\in \lbrack 0,T].
\end{align}%
We notice that due to the properties of the stochastic integral and the
assumption \eqref{G}, we have 
\begin{equation*}
\mathbf{h}_{n}(t)\rightharpoonup \mathbf{h}(t)=\xi _{0}^{2}(t)\mathbf{y}%
(t)-\int_{0}^{t}\xi _{0}^{2}(s)\mathbf{G}^{\ast }(t)\,d{\mathcal{W}}%
_{s}\quad \mbox{ weakly in
}\ L_{2}(\Omega \times (0,T);H^{1}({\mathcal{O}})).
\end{equation*}%
Now, we pass to the limit in the equation \eqref{CC1} in the distributional
sense. Namely multiplying the equation \eqref{11:33} by the test funtion $%
\theta (t)\eta (\omega )$, with $\theta \in C^{\infty }([0,T])$ with compact
support and $\eta \in L_{2}(\Omega )$, and passing to the limit, we derive 
\begin{align*}
\mathbb{E}\int_{0}^{T}\left( \mathbf{h}(t),\boldsymbol{\varphi }\right)
\theta ^{\prime }(t)\eta & =-\mathbb{E}\int_{0}^{T}\xi _{0}^{2}(t)[-\nu
\left( \mathbf{y}(t),\boldsymbol{\varphi }\right) _{V}\,+\nu \int_{\Gamma
}b(t)(\boldsymbol{\varphi }\cdot {\bm{\tau }})\,d\mathbf{\gamma }\, \\
& -\left( B^{\ast }(t),\boldsymbol{\varphi }\right) -2f_{0}(t)\left( \mathbf{%
y}(t),\boldsymbol{\varphi }\right) ]\theta \eta \,dt.
\end{align*}%
Therefore $\frac{\partial \mathbf{h}}{\partial t}\in L_{2}(\Omega \times
(0,T);\left( H^{1}(\mathcal{O})\right) ^{\ast })$. Since $\mathbf{h}\in
L_{2}(\Omega \times (0,T);H^{1}(\mathcal{O}))$, we infer that $\mathbf{h}\in
L_{2}(\Omega ;C([0,T];L_{2}(\mathcal{O}))$ by the Aubin-Lions embeeding
result \cite{aub}, \cite{tem}. Taking into account the continuity property
of the stochastic integral, we conclude that $\xi _{0}^{2}\mathbf{y}\in
L_{2}(\Omega ;C([0,T];L_{2}(\mathcal{O})).$ In addition 
\begin{equation*}
\xi _{0}^{2}\mathbf{y}_{n}\rightharpoonup \xi _{0}^{2}\mathbf{y}\qquad \text{%
in}\quad C_{\omega }([0,T],L_{2}(\Omega )\times L_{2}(\mathcal{O})),
\end{equation*}%
where the index $\omega $ means that we are considering $L_{2}(\Omega
)\times L_{2}(\mathcal{O})$ endowed with the weak topology. Hence, we have 
\begin{equation}
\mathbb{E}\left[ \left( \xi _{0}^{2}(t)\mathbf{y}_{n}(t),\boldsymbol{\varphi 
}\right) \eta \right] \rightarrow \mathbb{E}\left[ \left( \xi _{0}^{2}(t)%
\mathbf{y}(t),\boldsymbol{\varphi }\right) \eta \right] ,\qquad \forall t\in
\lbrack 0,T].  \label{ER}
\end{equation}

Now, we multiply the equation \eqref{CC1} by an arbitrary $\eta \in
L_{2}(\Omega )$ and take the expectation, we derive%
\begin{eqnarray*}
&&\mathbb{E\ }\eta \{\left( \xi _{0}^{2}(t)\mathbf{y}_{n}(t),\boldsymbol{%
\varphi }\right) -\left( \mathbf{y}_{n,0},\boldsymbol{\varphi }\right) \} \\
&&\qquad \qquad =\mathbb{E\ }\eta \left\{ \int_{0}^{t}\xi _{0}^{2}[-\nu
\left( \mathbf{y}_{n},\boldsymbol{\varphi }\right) _{V}\,+\nu \int_{\Gamma
}b(\boldsymbol{\varphi }\cdot {\bm{\tau }})\,d\mathbf{\gamma }-(\left( 
\mathbf{y}_{n}\cdot \nabla )\mathbf{y}_{n},\boldsymbol{\varphi }\right)
\,\right. \\
&&\qquad \qquad \left. -2f_{0}\left( \mathbf{y}_{n},\boldsymbol{\varphi }%
\right) ]\,dt+\int_{0}^{t}\xi _{0}^{2}\left( \mathbf{G}(s,\mathbf{y}_{n}),%
\boldsymbol{\varphi }\right) \,d{\mathcal{W}}_{s}\right\} .
\end{eqnarray*}%
Applying (\ref{c1})-(\ref{c01}) and (\ref{ER}), we pass to the limit $%
n\rightarrow \infty $ in this equality and deduce%
\begin{eqnarray*}
&&\mathbb{E\ }\eta \left\{ \left( \xi _{0}^{2}(t)\mathbf{y}(t),\boldsymbol{%
\varphi }\right) -\left( \mathbf{y}_{0},\boldsymbol{\varphi }\right)
\right\} =\mathbb{E\ }\eta \left\{ \int_{0}^{t}\xi _{0}^{2}[-\nu \left( 
\mathbf{y},\boldsymbol{\varphi }\right) _{V}\,+\nu \int_{\Gamma }b(%
\boldsymbol{\varphi }\cdot {\bm{\tau }})\,d\mathbf{\gamma }\right. \\
&&\left. \quad \quad -\left( B^{\ast },\boldsymbol{\varphi }\right)
-2f_{0}\left( \mathbf{y},\boldsymbol{\varphi }\right) ]\,dt+\int_{0}^{t}\xi
_{0}^{2}\left( \mathbf{G}^{\ast }(s),\boldsymbol{\varphi }\right) \,d{%
\mathcal{W}}_{s}\right\} .
\end{eqnarray*}%
Since $\eta \in L_{2}(\Omega )$ is arbitrary, the following equation holds 
\begin{align}
& \left( \xi _{0}^{2}(t)\mathbf{y}(t),\boldsymbol{\varphi }\right) -\left( 
\mathbf{y}_{0},\boldsymbol{\varphi }\right) =\left\{ \int_{0}^{t}\xi
_{0}^{2}[-\nu \left( \mathbf{y},\boldsymbol{\varphi }\right) _{V}\,+\nu
\int_{\Gamma }b(\boldsymbol{\varphi }\cdot {\bm{\tau }})\,d\mathbf{\gamma }%
\right.  \notag \\
& \left. \quad \quad -\left( B^{\ast },\boldsymbol{\varphi }\right)
-f_{0}\left( \mathbf{y},\boldsymbol{\varphi }\right) ]\,dt+\int_{0}^{t}\xi
_{0}^{2}\left( \mathbf{G}^{\ast },\boldsymbol{\varphi }\right) \,d{\mathcal{W%
}}_{s}\right\}  \label{555}
\end{align}%
for any $t\in \lbrack 0,T]$ and $P$-a.e. in $\Omega ,$ that is%
\begin{align*}
d\left( \xi _{0}^{2}\left( \mathbf{y},\boldsymbol{\varphi }\right) \right) &
=\xi _{0}^{2}[-\nu \left( \mathbf{y},\boldsymbol{\varphi }\right) _{V}\,+\nu
\int_{\Gamma }b(\boldsymbol{\varphi }\cdot {\bm{\tau }})\,d\mathbf{\gamma }%
-\left( B^{\ast },\boldsymbol{\varphi }\right) \, \\
& -2f_{0}\left( \mathbf{y},\boldsymbol{\varphi }\right) ]\,dt+\xi
_{0}^{2}\left( \mathbf{G}^{\ast },\boldsymbol{\varphi }\right) \,d{\mathcal{W%
}}_{t}\text{\quad and\quad }\mathbf{y}(0)=\mathbf{y}_{0}.
\end{align*}

Moreover if we use It\^o's formula $\ $%
\begin{equation*}
\ d\left( \mathbf{y},\boldsymbol{\varphi }\right) =d\left[ \xi _{0}^{-2}\xi
_{0}^{2}\left( \mathbf{y},\boldsymbol{\varphi }\right) \right] =\xi
_{0}^{2}\left( \mathbf{y},\boldsymbol{\varphi }\right) d\left( \xi
_{0}^{-2}\right) +\xi _{0}^{-2}d\left[ \xi _{0}^{2}\left( \mathbf{y},%
\boldsymbol{\varphi }\right) \right] ,
\end{equation*}%
we derive that the limit function $\mathbf{y}$ in the form $\mathbf{y}=%
\mathbf{u}+\mathbf{a}$ with 
\begin{equation*}
\mathbf{u\in }L_{\infty }(0,T;H)\cap L_{2}(0,T;V),\quad P\text{-a.e. in }%
\Omega ,\quad \text{a.e. on }(0,T),
\end{equation*}%
satisfies P-a.e. in $\Omega $\ the stochastic differential equation 
\begin{align}
d\left( \mathbf{y},\boldsymbol{\varphi }\right) & =[-\nu \left( \mathbf{y},%
\boldsymbol{\varphi }\right) _{V}+\nu \int_{\Gamma }b(\boldsymbol{\varphi }%
\cdot {\bm{\tau }})\,d\mathbf{\gamma }-(B^{\ast }(t),\boldsymbol{\varphi }%
)\,]dt  \notag \\
& +\left( \mathbf{G}^{\ast }(t),\boldsymbol{\varphi }\right) \,d{\mathcal{W}}%
_{t},\qquad \forall t\in \lbrack 0,T],\quad \forall \boldsymbol{\varphi }\in
V,  \label{y12}
\end{align}%
and $\mathbf{y}(0)=\mathbf{y}_{0}.$

\textit{Step 3.} \textit{Deduction of strong convergence as }$n\rightarrow
\infty .$ In order to prove that the limit process $\mathbf{y}$ satisfies
the equation \eqref{res1}, we adapt the methods in \cite{B99}. Writing $%
\mathbf{y}=\mathbf{u}+\mathbf{a}$, $\mathbf{y}_{n}=\mathbf{u}_{n}+\mathbf{a}$
and taking the difference of the equations (\ref{y1}) and (\ref{y12}) with $%
\boldsymbol{\varphi }=\mathbf{e}_{i}\in V_{n}$, $i=1,...,n$, we deduce 
\begin{eqnarray}
d\left( P_{n}\mathbf{u}-\mathbf{u}_{n},\mathbf{e}_{i}\right) &=&\left[ -\nu
\left( P_{n}\mathbf{u}-\mathbf{u}_{n},\mathbf{e}_{i}\right) _{V}+\left( (%
\mathbf{y}_{n}\cdot \nabla )\mathbf{y}_{n}-B^{\ast }(t),\mathbf{e}%
_{i}\right) \right] \,dt  \notag \\
&-&\left( \mathbf{G}(t,\mathbf{y}_{n})-\mathbf{G}^{\ast }(t),\mathbf{e}%
_{i}\right) \,d{\mathcal{W}}_{t},\quad i=1,...,n.  \label{diff1}
\end{eqnarray}%
Then the It\^{o}'s formula yields 
\begin{align*}
d(P_{n}\mathbf{u}-\mathbf{u}_{n},\mathbf{e}_{i})^{2}& =2\left( P_{n}\mathbf{u%
}-\mathbf{u}_{n},\mathbf{e}_{i}\right) \\
& \times \left[ -\nu \left( P_{n}\mathbf{u}-\mathbf{u}_{n},\mathbf{e}%
_{i}\right) _{V}+\left( (\mathbf{y}_{n}\cdot \nabla )\mathbf{y}_{n}-B^{\ast
}(t),\mathbf{e}_{i}\right) \right] \,dt\vspace{2mm} \\
& -2\left( P_{n}\mathbf{u}-\mathbf{u}_{n},\mathbf{e}_{i}\right) \left( 
\mathbf{G}(t,\mathbf{y}_{n})-\mathbf{G}^{\ast }(t),\mathbf{e}_{i}\right) \,d{%
\mathcal{W}}_{t} \\
& +|\left( \mathbf{G}(t,\mathbf{y}_{n})-\mathbf{G}^{\ast }(t),\mathbf{e}%
_{i}\right) |^{2}\,dt.
\end{align*}%
Summing over $i=1,\dots ,n,$ we derive%
\begin{align}
d\left( ||P_{n}\mathbf{u}-\mathbf{u}_{n}||_{2}^{2}\right) & +2\nu ||P_{n}%
\mathbf{u}-\mathbf{u}_{n}||_{V}^{2}dt=2((\mathbf{y}_{n}\cdot \nabla )\mathbf{%
y}_{n}-B^{\ast }(t),P_{n}\mathbf{u}-\mathbf{u}_{n})\,dt  \notag \\
& +\sum_{i=1}^{n}|\left( \mathbf{G}(t,\mathbf{y}_{n})-\mathbf{G}^{\ast }(t),%
\mathbf{e}_{i}\right) |^{2}\,dt  \notag \\
& -2\left( \mathbf{G}(t,\mathbf{y}_{n})-\mathbf{G}^{\ast }(t),P_{n}\mathbf{u}%
-\mathbf{u}_{n}\right) \,d{\mathcal{W}}_{t}.  \label{y13}
\end{align}

Standard computations give 
\begin{align*}
(\mathbf{y}_{n}\cdot \nabla )\mathbf{y}_{n}-B^{\ast }(t)& =\left\{ -((%
\mathbf{u}_{n}+\mathbf{a})\cdot \nabla )(P_{n}\mathbf{u}-\mathbf{u}%
_{n})-((P_{n}\mathbf{u}-\mathbf{u}_{n})\cdot \nabla )(P_{n}\mathbf{u}+%
\mathbf{a})\right\} \\
& +((P_{n}\mathbf{u}-\mathbf{u})\cdot \nabla )(P_{n}\mathbf{u}+\mathbf{a})+((%
\mathbf{u}+\mathbf{a})\cdot \nabla )(P_{n}\mathbf{u}-\mathbf{u}) \\
& +\left\{ (\mathbf{y}\cdot \nabla )\mathbf{y}-B^{\ast }(t)\right\} =\left\{
A_{0,1}+A_{0,2}\right\} +A_{1}+A_{2}+A_{3}.
\end{align*}%
In addition, using (\ref{LI}), (\ref{TT}), (\ref{cal}) and Theorem 4.47, p.
210, of \cite{du}, we show the existence of a constant $C_{2}$, verifying
the relation 
\begin{align}
I_{0}& :=|\left( \left\{ A_{0,1}+A_{0,2}\right\} ,P_{n}\mathbf{u}-\mathbf{u}%
_{n}\right) |  \notag \\
& \leqslant \Bigl|\int_{\Gamma }a((P_{n}\mathbf{u}-\mathbf{u}_{n})\cdot %
\bm{\tau })^{2}\,d\mathbf{\gamma }\Bigr|+\left\vert \left( ((P_{n}\mathbf{u}-%
\mathbf{u}_{n})\cdot \nabla )(P_{n}\mathbf{u}+\mathbf{a}),P_{n}\mathbf{u}-%
\mathbf{u}_{n}\right) \right\vert  \notag \\
& \leqslant \Vert a\Vert _{L_{\infty }(\Gamma )}\Vert P_{n}\mathbf{u}-%
\mathbf{u}_{n}\Vert _{L_{2}(\Gamma )}^{2}+\left\Vert P_{n}\mathbf{u}+\mathbf{%
\ a}\right\Vert _{V}\left\Vert P_{n}\mathbf{u}-\mathbf{u}_{n}\right\Vert
_{4}^{2}  \notag \\
& \leqslant (\Vert a\Vert _{L_{\infty }(\Gamma )}+\left\Vert \mathbf{a}%
\right\Vert _{H^{1}}+\left\Vert P_{n}\mathbf{u}\right\Vert _{V})\left\Vert
P_{n}\mathbf{u}-\mathbf{u}_{n}\right\Vert _{2}\Vert P_{n}\mathbf{u}-\mathbf{u%
}_{n}\Vert _{V}  \notag \\
& \leqslant C_{2}(||(a,b)||_{\mathcal{H}_{p}(\Gamma )}^{2}+\left\Vert 
\mathbf{u}\right\Vert _{V}^{2})\Vert P_{n}\mathbf{u}-\mathbf{u}_{n}\Vert
_{2}^{2}+\nu \left\Vert P_{n}\mathbf{u}-\mathbf{u}_{n}\right\Vert _{V}^{2}.
\label{BB1}
\end{align}%
On the other hand, H\"{o}lder's inequality gives%
\begin{align}
I_{1}& :=|\left( A_{1},P_{n}\mathbf{u}-\mathbf{u}_{n}\right) |\leqslant
C||P_{n}\mathbf{u}-\mathbf{u}||_{4}||\nabla (P_{n}\mathbf{u}+\mathbf{a}%
)||_{2}||P_{n}\mathbf{u}-\mathbf{u}_{n}||_{4}  \notag \\
& \leqslant C||P_{n}\mathbf{u}-\mathbf{u}||_{4}(||\mathbf{u}||_{V}+||\mathbf{%
a}||_{H^{1}})(||P_{n}\mathbf{u}||_{4}+||\mathbf{u}_{n}||_{4})  \label{BB22}
\end{align}%
and%
\begin{align}
I_{2}& :=|\left( A_{2},P_{n}\mathbf{u}-\mathbf{u}_{n}\right) |\leqslant C||%
\mathbf{u}+\mathbf{a}||_{4}||\nabla (P_{n}\mathbf{u}-\mathbf{u})||_{2}||P_{n}%
\mathbf{u}-\mathbf{u}_{n}||_{4}  \notag \\
& \leqslant C||P_{n}\mathbf{u}-\mathbf{u}||_{V}(||\mathbf{u}||_{4}+||\mathbf{%
a}||_{4})(||P_{n}\mathbf{u}||_{4}\mathbf{+||u}_{n}||_{4}).  \label{BB22_2}
\end{align}%
The last term $A_{3}$ will be considered later on.

Now, by denoting 
\begin{equation}
\mathbf{G}_{n}=\mathbf{G}(t,\mathbf{y}_{n}),\quad \mathbf{G}=\mathbf{G}(t,%
\mathbf{y}),\quad \mathbf{G}^{\ast }=\mathbf{G}^{\ast }(t),  \label{GS_NS}
\end{equation}%
we have%
\begin{equation*}
\sum_{i=1}^{n}|\left( \mathbf{G}(t,\mathbf{y}_{n})-\mathbf{G}^{\ast }(t),%
\mathbf{e}_{i}\right) |^{2}=\sum_{i=1}^{n}|(\mathbf{G}_{n}-\mathbf{G}^{\ast
},\mathbf{e}_{i})|^{2}=\Vert P_{n}\mathbf{G}_{n}-P_{n}\mathbf{G}^{\ast
}\Vert _{2}^{2}.
\end{equation*}%
The standard relation $x^{2}=(x-y)^{2}-y^{2}+2xy\;$ allows to write 
\begin{align*}
\Vert P_{n}\mathbf{G}_{n}-P_{n}\mathbf{G}^{\ast }\Vert _{2}^{2}& =\Vert P_{n}%
\mathbf{G}_{n}-P_{n}\mathbf{G}\Vert _{2}^{2}-\Vert P_{n}\mathbf{G}-P_{n}%
\mathbf{G}^{\ast }\Vert _{2}^{2} \\
& -2(P_{n}\mathbf{G}_{n}-P_{n}\mathbf{G}^{\ast },P_{n}\mathbf{G}-P_{n}%
\mathbf{G}^{\ast }).
\end{align*}%
From (\ref{G}) and (\ref{2})$_{1}$, we have 
\begin{equation*}
\Vert P_{n}\mathbf{G}_{n}-P_{n}\mathbf{G}\Vert _{2}^{2}\leqslant \Vert 
\mathbf{G}_{n}-\mathbf{G}\Vert _{2}^{2}\leqslant K\left\Vert \mathbf{u}_{n}-%
\mathbf{u}\right\Vert _{2}^{2},
\end{equation*}%
then for the fixed constant $C_{3}=2K$ it follows that 
\begin{align}
\sum_{i=1}^{n}|(\mathbf{G}(t,\mathbf{y}_{n})& -\mathbf{G}^{\ast }(t),\mathbf{%
e}_{i})|^{2}=\Vert P_{n}\mathbf{G}_{n}-P_{n}\mathbf{G}^{\ast }\Vert _{2}^{2}
\notag \\
& \leqslant K\left\Vert \mathbf{u}_{n}-\mathbf{u}\right\Vert _{2}^{2}-\Vert
P_{n}\mathbf{\ G}-P_{n}\mathbf{G}^{\ast }\Vert _{2}^{2}  \notag \\
& +2(P_{n}\mathbf{G}_{n}-P_{n}\mathbf{G}^{\ast },P_{n}\mathbf{G}-P_{n}%
\mathbf{G}^{\ast })  \notag \\
& \leqslant C_{3}\left\Vert \mathbf{u}_{n}-P_{n}\mathbf{u}\right\Vert
_{2}^{2}+C\left\Vert P_{n}\mathbf{u}-\mathbf{u}\right\Vert _{2}^{2}-\Vert
P_{n}\mathbf{G}-P_{n}\mathbf{G}^{\ast }\Vert _{2}^{2}  \notag \\
& +2(P_{n}\mathbf{G}_{n}-P_{n}\mathbf{G}^{\ast },P_{n}\mathbf{G}-P_{n}%
\mathbf{G}^{\ast }).  \label{GG1}
\end{align}%
The positive constants $C_{2}$ and $C_{3}$ in (\ref{BB1}) and (\ref{GG1})
are independent of $n.$

We notice that with the help of the convergence results (\ref{c02}), (\ref%
{c1})-(\ref{c01}), and performing a suitable limit trfansition in \ the
equation (\ref{y13}), as $n\rightarrow \infty ,$ we can verify that all
terms on the right-hand side of the equality (\ref{y13}) containing $P_{n}%
\mathbf{u}-\mathbf{u}$ will vanish; 
% according to relations (\ref{c02}) and (\ref{BB22})-(\ref{BB22_2}), (%\ref{GG1})
however, terms that contain $P_{n}\mathbf{u}-\mathbf{u}_{n}$ \ will remain.
Fortunately, these terms can be eliminated by introducing the auxiliary
function 
\begin{equation}
\widetilde{\xi }(t)=e^{-\int_{0}^{t}\widetilde{f}(s)\ ds}  \label{kksi}
\end{equation}%
with $\widetilde{f}(t)=C_{3}+\max (3C_{0},C_{2})(1+||(a,b)||_{\mathcal{H}%
_{p}(\Gamma )}^{2}+\left\Vert \mathbf{u}\right\Vert _{V}^{2}).$

Now, by applying It\^{o}'s formula to the equality (\ref{y13}) and using the
definition (\ref{kksi}) of $\widetilde{\xi }$, we obtain%
\begin{align*}
d(\widetilde{\xi }(t)||P_{n}\mathbf{u}& -\mathbf{u}_{n}||_{2}^{2})+2\nu 
\widetilde{\xi }(t)||P_{n}\mathbf{u}-\mathbf{u}_{n}||_{V}^{2}\,dt \\
& \leqslant 2\widetilde{\xi }(t)((\mathbf{y}_{n}\cdot \nabla )\mathbf{y}%
_{n}-B^{\ast }(t),P_{n}\mathbf{u}-\mathbf{u}_{n})\,dt \\
& +\widetilde{\xi }(t)\sum_{i=1}^{n}|\left( \mathbf{G}(t,\mathbf{y}_{n})-%
\mathbf{G}^{\ast }(t),\mathbf{e}_{i}\right) |^{2}\,dt \\
& -2\widetilde{\xi }(t)\left( \mathbf{G}(t,\mathbf{y}_{n})-\mathbf{G}^{\ast
}(t),P_{n}\mathbf{u}-\mathbf{u}_{n}\right) \,d{\mathcal{W}}_{t}-C_{3}%
\widetilde{\xi }(t)||P_{n}\mathbf{u}-\mathbf{u}_{n}||_{2}^{2}\,dt \\
& -C_{2}\widetilde{\xi }(t)(||(a,b)||_{\mathcal{H}_{p}(\Gamma
)}^{2}+\left\Vert \mathbf{u}\right\Vert _{V}^{2})||P_{n}\mathbf{u}-\mathbf{u}%
_{n}||_{2}^{2}\,dt.
\end{align*}%
Writing this equation in the integral form, taking the expectation, and
applying the estimates (\ref{BB1}), (\ref{GG1}), we deduce that 
\begin{align*}
\mathbb{E}(\widetilde{\xi }(t)||P_{n}\mathbf{u}(t)& -\mathbf{u}%
_{n}(t)||_{2}^{2})+\mathbb{E}\int_{0}^{t}\widetilde{\xi }(s)\Vert P_{n}%
\mathbf{G}-P_{n}\mathbf{G}^{\ast }\Vert _{2}^{2}ds \\
& +\nu \mathbb{\ E}\int_{0}^{t}\widetilde{\xi }(s)||P_{n}\mathbf{u}-\mathbf{u%
}_{n}||_{V}^{2}\,ds\leqslant 2\mathbb{E}\int_{0}^{t}\widetilde{\xi }%
(s)I_{1}\,ds \\
& +2\mathbb{E}\int_{0}^{t}~\widetilde{\xi }(s)I_{2}\,ds+2\mathbb{E}%
\int_{0}^{t}\widetilde{\xi }(s)\left( A_{3},P_{n}\mathbf{u}-\mathbf{u}%
_{n}\right) \,ds \\
& +C\mathbb{E}\int_{0}^{t}\widetilde{\xi }(s)\left\Vert P_{n}\mathbf{u}-%
\mathbf{\ u}\right\Vert _{2}^{2}\,ds \\
& +2\mathbb{E}\int_{0}^{t}\widetilde{\xi }(s)(P_{n}\mathbf{G}_{n}-P_{n}%
\mathbf{G}^{\ast },P_{n}\mathbf{G}-P_{n}\mathbf{G}^{\ast })\,ds \\
& =J_{1}+J_{2}+J_{3}+J_{4}+J_{5}\qquad \text{for \ }t\in (0,T).
\end{align*}%
Next, we will show that the right-hand side of this inequality tends to zero
as $n\rightarrow \infty $.

Considering the estimate (\ref{BB22}) and using $\widetilde{\xi }\leqslant
\xi _{0}^{3}$ on $(0,T)$, then we deduce that%
\begin{eqnarray*}
J_{1} &\leqslant &C\left( \mathbb{E}\int_{0}^{T}\xi _{0}^{3}||P_{n}\mathbf{u}%
-\mathbf{u}||_{4}^{2}(||\mathbf{u}||_{V}+||\mathbf{a}||_{H^{1}})\ ds\right)
^{1/2} \\
&\times &\left( \mathbb{E}\int_{0}^{T}\xi _{0}^{3}(||\mathbf{u}||_{V}+||%
\mathbf{a}||_{H^{1}})(||P_{n}\mathbf{u}||_{4}^{2}+||\mathbf{u}%
_{n}||_{4}^{2})\ ds\right) ^{1/2}.
\end{eqnarray*}%
Using (\ref{LI}) for $q=4$, we have 
\begin{eqnarray*}
\mathbb{E}\int_{0}^{T}\xi _{0}^{3}||P_{n}\mathbf{u} &-&\mathbf{u}%
||_{4}^{2}(||\mathbf{u}||_{V}+||\mathbf{a}||_{H^{1}})\ ds \\
&\leqslant &(\mathbb{E}\sup_{s\in \lbrack 0,t]}\xi _{0}^{2}||P_{n}\mathbf{u}-%
\mathbf{u}||_{2}^{2}\int_{0}^{T}\xi _{0}^{2}(||\mathbf{u}||_{V}^{2}+||%
\mathbf{a}||_{H^{1}}^{2})\ ds)^{1/2} \\
&\times &(\mathbb{E}\int_{0}^{T}\xi _{0}^{2}||P_{n}\mathbf{u}-\mathbf{u}%
||_{V}^{2}\ ds)^{1/2}\leqslant C(\mathbb{E}\int_{0}^{T}\xi _{0}^{2}||P_{n}%
\mathbf{u}-\mathbf{u}||_{V}^{2}\ ds)^{1/2}
\end{eqnarray*}%
by the estimates (\ref{21:47})-(\ref{uny}). Applying similar calculations we
can show that there exists a constant $C$, such that 
\begin{equation*}
\mathbb{E}\int_{0}^{T}\xi _{0}^{3}(||\mathbf{u}||_{V}+||\mathbf{a}%
||_{H^{1}})(||P_{n}\mathbf{u}||_{4}^{2}+||\mathbf{u}_{n}||_{4}^{2})\
ds\leqslant C,
\end{equation*}%
that is%
\begin{equation*}
J_{1}\leqslant C\left( \mathbb{E}\int_{0}^{T}\xi _{0}^{2}||P_{n}\mathbf{u}-%
\mathbf{u}||_{V}^{2}\ ds\right) ^{1/4}.
\end{equation*}%
For the term $J_{2},$\ \ using the estimate (\ref{BB22_2}), we can show that%
\begin{equation*}
J_{2}\leqslant C\left( \mathbb{E}\int_{0}^{T}\xi _{0}^{2}||P_{n}\mathbf{u}-%
\mathbf{u}||_{V}^{2}\ ds\right) ^{1/2}.
\end{equation*}%
Therefore we get that the terms $J_{i},$ $i=1,2,$\ converge to zero as $%
n\rightarrow \infty $ by (\ref{c02Y}).

The convergences of (\ref{c1}) and (\ref{c02Y}) show that 
\begin{equation*}
\xi _{0}\left( P_{n}\mathbf{u}-\mathbf{u}_{n}\right) \rightharpoonup 0\qquad 
\text{ weakly in }L_{2}(\Omega \times (0,T),V)\qquad \text{as }n\rightarrow
\infty .
\end{equation*}%
The operator $\xi _{0}^{2}A_{3}=\xi _{0}^{2}\left( (\mathbf{y}\cdot \nabla )%
\mathbf{y}-B^{\ast }\right) $ belongs to $L_{2}(\Omega \times
(0,T);V^{\prime })$\ by (\ref{l01}) and (\ref{c01}), thus 
\begin{equation*}
J_{3}=2\mathbb{E}\int_{0}^{T}\widetilde{\xi }(s)((\mathbf{y}\cdot \nabla )%
\mathbf{y})-B^{\ast },P_{n}\mathbf{u}-\mathbf{u}_{n})\,ds\rightarrow 0\qquad 
\text{as }n\rightarrow \infty .
\end{equation*}

Due to (\ref{c02Y}), we have%
\begin{equation*}
J_{4}=C\mathbb{E}\int_{0}^{T}\widetilde{\xi }(s)\left\Vert P_{n}\mathbf{u}-%
\mathbf{u}\right\Vert _{2}^{2}ds\rightarrow 0.
\end{equation*}

Due to the convergence results (\ref{c02}), (\ref{c1}), (\ref{c02Y}), (\ref%
{c01}) and (\ref{GS_NS}), we obtain 
\begin{eqnarray}
\xi _{0}P_{n}\left( \mathbf{G}_{n}-\mathbf{G}^{\ast }\right)
&\rightharpoonup &\mathbf{0}\text{\qquad\ weakly in }L_{2}(\Omega \times
(0,T),H^{m}),  \notag \\
\xi _{0}P_{n}(\mathbf{G}-\mathbf{G}^{\ast }) &\rightarrow &\mathbf{G}-%
\mathbf{G}^{\ast }\text{\qquad strongly in }L_{2}(\Omega \times (0,T),H^{m}),
\label{11L}
\end{eqnarray}%
that implies 
\begin{equation*}
J_{5}=2\mathbb{E}\int_{0}^{T}\widetilde{\xi }(s)(P_{n}\mathbf{G}_{n}-P_{n}%
\mathbf{G}^{\ast },P_{n}(\mathbf{G}-\mathbf{G}^{\ast }))\,ds\rightarrow 0,%
\text{ as }n\rightarrow \infty .
\end{equation*}%
After combining all the convergence results, we obtain the following strong
convergences\ 
\begin{equation*}
\lim_{n\rightarrow \infty }\mathbb{E}\left( \widetilde{\xi }(t)||P_{n}%
\mathbf{u}(t)-\mathbf{u}_{n}(t)||_{2}^{2}\right) =0,\text{\qquad }%
\lim_{n\rightarrow \infty }\mathbb{E}\int_{0}^{t}\widetilde{\xi }(s)||P_{n}%
\mathbf{u}-\mathbf{u}_{n}||_{V}^{2}\,ds=0
\end{equation*}%
for $t\in (0,T)$, which combined with (\ref{c02Y}), imply 
\begin{equation}
\lim_{n\rightarrow \infty }\mathbb{E}\left( \widetilde{\xi }(t)||\mathbf{u}%
_{n}(t)-\mathbf{u}(t)||_{2}^{2}\right) =0,\text{\qquad }\lim_{n\rightarrow
\infty }\mathbb{E}\int_{0}^{t}\widetilde{\xi }(s)||\mathbf{u}_{n}-\mathbf{u}%
||_{V}^{2}\,ds=0.  \label{gg}
\end{equation}

In addition, considering (\ref{G}), we conclude%
\begin{equation*}
\mathbb{E}\int_{0}^{t}\widetilde{\xi }(s)\Vert \mathbf{G}(s,\mathbf{y})-%
\mathbf{G}^{\ast }(s)\Vert _{2}^{2}ds=0.
\end{equation*}%
Since $\widetilde{\xi }$ is strictly positive, we infer that 
\begin{equation}
\mathbf{G}(t,\mathbf{y})=\mathbf{G}^{\ast }(t)\qquad \text{a. e. in }\Omega
\times (0,T).  \label{099}
\end{equation}%
From (\ref{c01}) and (\ref{gg}), it follows that $\widetilde{\xi }(t)(%
\mathbf{y}\cdot \nabla )\mathbf{y}=\widetilde{\xi }(t)\mathbf{B}^{\ast }(t)\ 
$ a. e. in $\Omega \times (0,T),$ that implies 
\begin{equation}
(\mathbf{y}\cdot \nabla )\mathbf{y}=\mathbf{B}^{\ast }(t)\qquad \text{a. e.
in }\Omega \times (0,T).  \label{RC2YY}
\end{equation}%
Considering the identities (\ref{099}), ( \ref{RC2YY}), the equation (\ref{555})
reads 
\begin{align*}
\left( \mathbf{y}(t),\boldsymbol{\varphi }\right) -\left( \mathbf{y}_{0},%
\boldsymbol{\varphi }\right) & =\int_{0}^{t}\left[ -\nu \left( \mathbf{y},%
\boldsymbol{\varphi }\right) _{V}+\nu \int_{\Gamma }b(\boldsymbol{\varphi }%
\cdot {\bm{\tau }})\,d\mathbf{\gamma }-(\mathbf{y}\cdot \nabla )\mathbf{y},%
\boldsymbol{\varphi })\right] \,ds \\
& +\int_{0}^{t}\left( \mathbf{G}(s,\mathbf{y}),\boldsymbol{\varphi }\right)
\,d{\mathcal{W}}_{s},\quad \quad P\text{-a.e.in }\Omega ,\quad t\in (0,T).
\end{align*}

The uniqueness of the solution $\mathbf{y}$ follows from the stability
result established in the next theorem.
\end{proof}

Let us denote by $\widehat{\mathbf{\varphi }}=\mathbf{\varphi }_{1}-\mathbf{%
\varphi }_{2}$ the diference of two given functions $\mathbf{\varphi }_{1},%
\mathbf{\varphi }_{2}.$

\begin{theorem}
\label{Lips} Let us consider $\mathbf{y}_{1}=\mathbf{u}_{1}+\mathbf{a}_{1},$%
\ $\mathbf{y}_{2}=\mathbf{u}_{2}+\mathbf{a}_{2}$ with 
\begin{equation*}
\mathbf{u}_{1},\mathbf{u}_{2}\in C([0,T];H)\cap L_{4}(0,T;V),\quad P%
\text{-a.e.in }\Omega ,
\end{equation*}%
two solutions of \eqref{NSy}, satisfying the estimates $(\ref{21:47})$, $(%
\ref{uny})$ \ with two corresponding boundary conditions $a_{1},\;%
b_{1}$, $a_{2},$ $b_{2}$ and the initial
conditions 
\begin{equation*}
\mathbf{y}_{1,0}=\mathbf{u}_{1,0}+\mathbf{a}_{1}(0),\qquad \mathbf{y}_{2,0}=%
\mathbf{u}_{2,0}+\mathbf{a}_{2}(0).
\end{equation*}%
Then there exist a strictly positive function $f_{1}(t)\in L_{1}(0,T)$ $%
\quad P$-a.e.in $\Omega ,$\ depending only on the data, such that the
following estimate 
\begin{eqnarray}
\mathbb{E}\sup_{s\in \lbrack 0,t]}\xi _{1}^{2}(s)\left\Vert \widehat{\mathbf{%
y}}(s)\right\Vert _{2}^{2} &+&2\nu \int_{0}^{t}\xi _{1}^{2}(s)\Vert \widehat{%
\mathbf{y}}(s)\Vert _{V}^{2}\,ds  \notag \\
&\leqslant &C(\mathbb{E}\left\Vert \widehat{\mathbf{y}}_{0}\right\Vert
_{2}^{2}+\mathbb{E}\int_{0}^{t}\xi _{1}^{2}||(\widehat{a},\widehat{b})||_{%
\mathcal{H}_{p}(\Gamma )}^{2}\,ds)  \label{lip_H3}
\end{eqnarray}%
is valid with the function $\xi _{1}$ defined as 
\begin{equation}
\xi _{1}(t)=e^{-\int_{0}^{t}f_{1}(s)ds}\qquad \text{with }f_{1}\in
L_{1}(0,T)\quad P\text{-a.e. in }\Omega .  \label{aux}
\end{equation}%
$\hfill $
\end{theorem}

\begin{proof}
The proof follows the same reasoning as the proof of Theorem \ref{the_1}.
\end{proof}

\section{\protect\bigskip Solution to the control problem}

\label{sec3}\setcounter{equation}{0} This section studies the existence of
an optimal solution to the optimal control problem $(\mathcal{P}).$ We
intend to control the solution of the system (\ref{NSy}) by boundary values $%
(a,b)$, which belongs to the space $\mathcal{A}$ of admissible controls
defined as a {\it compact} subset of $L_{2}(\Omega \times (0,T);\mathcal{H}%
_{p}(\Gamma ))$ verifying an exponential integrability condition. More
precisely, we assume that there exists a constant $\lambda >0$ such that 
\begin{equation}
\mathbb{E}e^{4C_{0}\int_{0}^{T}||(a,b)||_{\mathcal{H}_{p}(\Gamma
)}^{2}ds}<\lambda ,\quad \forall (a,b)\in \mathcal{A}.  \label{ksi2}
\end{equation}

\begin{remark}
\label{R4.1}
We notice that given a control pair $(a,b)\in \mathcal{A}$, the
corresponding state $\mathbf{y}=\mathbf{u}+\mathbf{a}$ defined as the
solution of the state equation \eqref{res1} belongs to $L_{2}(\Omega \times
(0,T)\times \mathcal{O})$. Namely, considering the auxiliar function $\xi
_{0}$ introduced in \eqref{ksi}, and the estimates \eqref{21:47}, \eqref{uny}%
, H\"{o}lder's inequality gives 
\begin{equation*}
\mathbb{E(}\sup_{t\in \lbrack 0,T]}\Vert u(t)\Vert _{2}^{2})\leqslant (%
\mathbb{E}\sup_{t\in \lbrack 0,T]}\Vert \xi _{0}(t)u(t)\Vert _{2}^{4})^{%
\frac{1}{2}}\left( \mathbb{E}\left( \xi _{0}^{-4}(T)\right) \right) ^{\frac{1%
}{2}}<\infty .
\end{equation*}%
Therefore, the cost functional \eqref{cost} is well defined for every $%
(a,b)\in \mathcal{A}.$
\end{remark}

Now, we write one of the main result of the article, which establishes the
existence of a solution for the optimal control problem $(\mathcal{P})$.

\begin{teo}
\label{main_existence} Assume that $(a,b)$ and $\mathbf{y}_{0}$
verify the regularity \eqref{eq00sec12}, \eqref{0}, such that $(a,b)$
belongs to the space $\mathcal{A}$. Then there exists at least one solution
for the optimal control problem $(\mathcal{P}).$
\end{teo}

\begin{proof}
Let us consider a minimizing sequence 
\begin{equation*}
(a_{n},b_{n},\mathbf{y}_{n})\in \mathcal{A}\times L_{2}(\Omega ;L_{\infty
}(0,T;L_{2}({\mathcal{O}}))\cap L_{2}(0,T;H^{1}({\mathcal{O}})))
\end{equation*}%
of the cost functional $J$, namely 
\begin{equation*}
J(a_{n},b_{n},\mathbf{y}_{n})\rightarrow d=\inf (\mathcal{P})\quad \text{ as 
}n\rightarrow \infty ,
\end{equation*}%
and $\mathbf{y}_{n}$ is the weak solution of the system \eqref{NSy} for the
sequence $(a_{n},b_{n})\in \mathcal{A}$.%
\begin{align}
d\left( \mathbf{y}_{n}\,,\boldsymbol{\varphi }\right) & =\left[ -\nu \left( 
\mathbf{y}_{n},\boldsymbol{\varphi }\right) _{V}\,+\nu \int_{\Gamma }b_{n}(%
\boldsymbol{\varphi }\cdot {\bm{\tau }})\,d\mathbf{\gamma }-(\left( \mathbf{y%
}_{n}\cdot \nabla )\mathbf{y}_{n},\boldsymbol{\varphi }\right) \right] dt 
\notag \\
& +\left( \mathbf{G}(t,\mathbf{y}_{n}),\boldsymbol{\varphi }\right) \,d{%
\mathcal{W}}_{t},\quad \forall \boldsymbol{\varphi }\in V,\quad \quad P\text{%
-a.e. in }\Omega ,\quad \forall t\in (0,T),  \notag \\
\mathbf{u}_{n}(0)& =\mathbf{u}_{0}\in H,  \label{cp1}
\end{align}%
Due to the compactness of $\mathcal{A}$, there exists a subsequence, still
indexed by $n$, such that 
\begin{equation}
\left( a_{n},b_{n}\right) \rightarrow (a,b)\qquad {\text{ strongly in}}\quad
L_{2}(\Omega \times (0,T);\mathcal{H}_{p}(\Gamma )).  \label{y}
\end{equation}%
From Theorem 4.9., p. 94, of \cite{B11}, there exists a subsequence of $%
\left( a_{n},b_{n}\right) $, still denoted by $\left( a_{n},b_{n}\right) $,
and a function $h\in L_{2}(\Omega \times (0,T))$ such that 
\begin{equation}
\Vert (a,b)\Vert _{\mathcal{H}_{p}(\Gamma )}\leqslant h\text{,\quad }\Vert
(a_{n},b_{n})\Vert _{\mathcal{H}_{p}(\Gamma )}\leqslant h,\text{\quad }%
\forall n\in \mathbb{N},\text{\quad a. e. in }\Omega \times (0,T).
\label{fh}
\end{equation}%
Considering the function $h=h(t)$, let us introduce the following weight%
\begin{equation}
\xi _{h}(t)=e^{-C_{0}(t+\int_{0}^{t}h^{2}(s)ds)},\qquad P\text{-a.e. in }%
\Omega \text{.}  \label{ksi22}
\end{equation}

If we replace $a,b,\mathbf{a}$ by $a_{n},b_{n},\mathbf{a}_{n}$, respectively
in the\ relations \eqref{ha}, \eqref{cal}, \ then, taking into account the
estimates \eqref{21:47}, \eqref{uny}, we conclude that the sequence $\mathbf{%
u}_{n}=\mathbf{y}_{n}-\mathbf{a}_{n}$, $n\in \mathbb{N},$ satisfies the
estimates%
\begin{equation*}
\mathbb{E}\sup_{s\in \lbrack 0,t]}\xi _{h}^{2}(s)\left\Vert \mathbf{u}%
_{n}(s)\right\Vert _{2}^{2}+\nu \mathbb{E}\int_{0}^{t}\xi
_{h}^{2}(s)\left\Vert \mathbf{u}_{n}\right\Vert _{V}^{2}\,ds\leqslant C,
\end{equation*}%
\begin{equation}
\mathbb{E}\sup_{s\in \lbrack 0,t]}\xi _{h}^{4}(s)\left\Vert \mathbf{u}%
_{n}(s)\right\Vert _{2}^{4}+8\nu ^{2}\mathbb{E(}\int_{0}^{t}\xi
_{h}^{2}(s)\left\Vert \mathbf{u}_{n}(s)\right\Vert _{V}^{2}ds)^{2}\leqslant C
\label{lp11}
\end{equation}%
for any $t\in \lbrack 0,T]$, where the constants $C$ are independent of $n.$
Therefore there exists a subsequence, still indexed by $n$, such that 
\begin{eqnarray}
\xi _{h}\mathbf{u}_{n} &\rightharpoonup &\xi _{h}\mathbf{u}\qquad 
\mbox{ weakly in
}L_{s}(\Omega ;L_{2}(0,T;V)),  \notag \\
\xi _{h}\mathbf{u}_{n} &\rightharpoonup &\xi _{h}\mathbf{u}\qquad 
\mbox{ *-weakly in
}L_{s}(\Omega ,L_{\infty }(0,T;H))\quad \text{for }s=2\text{ and }4.
\label{c11}
\end{eqnarray}%
In addition, the following uniform estimate holds 
\begin{equation}
\left\Vert \xi _{h}^{2}\left( \mathbf{y}_{n}\cdot \nabla \right) \mathbf{y}%
_{n}\right\Vert _{L_{2}(\Omega \times (0,T);V^{\prime })}\leqslant C,\quad
\forall n\in \mathbb{N}.  \label{l}
\end{equation}%
Hence there exist operators $B^{\ast }$ and $G^{\ast }$ such that 
\begin{align}
\xi _{h}^{2}\left( \mathbf{y}_{n}\cdot \nabla \right) \mathbf{y}_{n}&
\rightharpoonup \xi _{h}^{2}B^{\ast }(t)\qquad \mbox{ weakly in
}\ L_{2}(\Omega \times (0,T);V^{\prime }),  \notag \\
\xi _{h}\mathbf{G}(t,\mathbf{y}_{n})& \rightharpoonup \xi _{h}\mathbf{G}%
^{\ast }(t)\qquad \mbox{ weakly in
}\ L_{2}(\Omega \times (0,T);H^{m}).  \label{c011}
\end{align}

Arguments already used in Step 2 of the proof of Theorem \ref{the_1} allow
to pass to the limit equation \eqref{cp1} in the distributional sense, as $%
n\rightarrow \infty $, to obtain 
\begin{align}
d\left( \mathbf{y},\boldsymbol{\varphi }\right) & =[-\nu \left( \mathbf{y},%
\boldsymbol{\varphi }\right) _{V}+\nu \int_{\Gamma }b(\boldsymbol{\varphi }%
\cdot {\bm{\tau }})\,d\mathbf{\gamma }-(B^{\ast }(t),\boldsymbol{\varphi }%
)\,]dt  \notag \\
& +\left( \mathbf{G}^{\ast }(t),\boldsymbol{\varphi }\right) \,d{\mathcal{W}}%
_{s},\quad \quad P\text{-a.e. in }\Omega ,\quad \forall t\in (0,T),  \notag
\\
\mathbf{u}(0)& =\mathbf{u}_{0}\in H,\qquad \qquad \qquad \qquad \forall 
\boldsymbol{\varphi }\in V.  \label{cp2}
\end{align}

Writing $\mathbf{y}=\mathbf{u}+\mathbf{a}$, $\mathbf{y}_{n}=\mathbf{u}_{n}+%
\mathbf{a}_{n}$ and doing the difference between (\ref{cp1}) and (\ref{cp2})
with $\boldsymbol{\varphi }=\mathbf{e}_{i}$, $i\in \mathbb{N}$, we deduce 
\begin{eqnarray}
d\left( \mathbf{u}-\mathbf{u}_{n},\mathbf{e}_{i}\right) &=&\Bigl[-\nu \left( 
\mathbf{u}-\mathbf{u}_{n},\mathbf{e}_{i}\right) _{V}+\nu \int_{\Gamma
}(b-b_{n})(\boldsymbol{\mathbf{e}_{i}}\cdot {\bm{\tau }})\,d\mathbf{\gamma }
\notag \\
&&-\left( \partial _{t}(\mathbf{a}-\mathbf{a}_{n}),\mathbf{e}_{i}\right)
-\nu \left( \mathbf{a}-\mathbf{a}_{n},\mathbf{e}_{i}\right) _{V}  \notag \\
&&+\left( (\mathbf{y}_{n}\cdot \nabla )\mathbf{y}_{n}-B^{\ast }(t),\mathbf{e}%
_{i}\right) \Bigr]\,dt  \notag \\
&&\vspace{2mm}-\left( \mathbf{G}(t,\mathbf{y}_{n})-\mathbf{G}^{\ast }(t),%
\mathbf{e}_{i}\right) \,d{\mathcal{W}}_{t},  \label{diff}
\end{eqnarray}%
which holds for any element of the basis $\left\{ \mathbf{e}_{i}\right\} .$

By applying It\^{o}'s formula, the equation (\ref{diff}) gives 
\begin{align*}
d(\mathbf{u}& -\mathbf{u}_{n},\mathbf{e}_{i})^{2}=2\left( \mathbf{u}-\mathbf{%
u}_{n},\mathbf{e}_{i}\right) \Bigl[-\nu \left( \mathbf{u}-\mathbf{u}_{n},%
\mathbf{e}_{i}\right) _{V}+\nu \int_{\Gamma }(b-b_{n})(\boldsymbol{\mathbf{e}%
_{i}}\cdot {\bm{\tau }})\,d\mathbf{\gamma } \\
& -\left( \partial _{t}(\mathbf{a}-\mathbf{a}_{n}),\mathbf{e}_{i}\right)
-\nu \left( \mathbf{a}-\mathbf{a}_{n},\mathbf{e}_{i}\right) _{V} \\
& +\left( (\mathbf{y}_{n}\cdot \nabla )\mathbf{y}_{n}-B^{\ast }(t),\mathbf{e}%
_{i}\right) \Bigr]\,dt\vspace{2mm} \\
& -2\left( \mathbf{u}-\mathbf{u}_{n},\mathbf{e}_{i}\right) \left( \mathbf{G}%
(t,\mathbf{y}_{n})-\mathbf{G}^{\ast }(t),\mathbf{e}_{i}\right) \,d{\mathcal{W%
}}_{t}+|\left( \mathbf{G}(t,\mathbf{y}_{n})-\mathbf{G}^{\ast }(t),\mathbf{e}%
_{i}\right) |^{2}\,dt.
\end{align*}%
Summing over the index $i\in \mathbb{N}$, \ we derive%
\begin{align}
d(||\mathbf{u}& -\mathbf{u}_{n}||_{2}^{2})+2\nu ||\mathbf{u}-\mathbf{u}%
_{n}||_{V}^{2}dt=2((\mathbf{y}_{n}\cdot \nabla )\mathbf{y}_{n}-B^{\ast }(t),%
\mathbf{u}-\mathbf{u}_{n})\,dt  \notag \\
& +2\nu \int_{\Gamma }(b-b_{n})((\mathbf{u}-\mathbf{u}_{n})\cdot {\bm{\tau }}%
)\,d\mathbf{\gamma }  \notag \\
& -\left( \partial _{t}(\mathbf{a}-\mathbf{a}_{n}),\mathbf{u}-\mathbf{u}%
_{n}\right) -\nu \left( \mathbf{a}-\mathbf{a}_{n},\mathbf{u}-\mathbf{u}%
_{n}\right) _{V}\Bigr]\,dt\vspace{2mm}  \notag \\
& +\sum_{i=1}^{\infty }|\left( \mathbf{G}(t,\mathbf{y}_{n})-\mathbf{G}^{\ast
}(t),\mathbf{e}_{i}\right) |^{2}\,dt  \notag \\
& -2\left( \mathbf{G}(t,\mathbf{y}_{n})-\mathbf{G}^{\ast }(t),\mathbf{u}-%
\mathbf{u}_{n}\right) \,d{\mathcal{W}}_{t}.  \label{y131}
\end{align}%
We write 
\begin{align*}
(\mathbf{y}_{n}\cdot \nabla )\mathbf{y}_{n}-B^{\ast }(t)& =\left\{ -((%
\mathbf{u}_{n}+\mathbf{a}_{n})\cdot \nabla )(\mathbf{u}-\mathbf{u}_{n})-((%
\mathbf{u}-\mathbf{u}_{n})\cdot \nabla )(\mathbf{u}+\mathbf{a})\right\} \\
& -((\mathbf{u}_{n}+\mathbf{a}_{n})\cdot \nabla )(\mathbf{a}_{n}-\mathbf{a}%
)+((\mathbf{a}_{n}-\mathbf{a})\cdot \nabla )(\mathbf{u}+\mathbf{a}) \\
& +\left\{ (\mathbf{y}\cdot \nabla )\mathbf{y}-B^{\ast }(t)\right\} \\
& =B_{0}+B_{1}+B_{2}+B_{3}.
\end{align*}%
With the help of (\ref{LI}), (\ref{TT}), (\ref{cal}) and \eqref{lp11}, we
deduce the following estimates 
\begin{align}
I_{0}& =|\left( B_{0},\mathbf{u}-\mathbf{u}_{n}\right) |\leqslant
|\int_{\Gamma }a_{n}((\mathbf{u}-\mathbf{u}_{n})\cdot \bm{\tau })^{2}\,d%
\mathbf{\gamma }|  \notag \\
& +|\left( ((\mathbf{u}-\mathbf{u}_{n})\cdot \nabla )(\mathbf{u}+\mathbf{a}),%
\mathbf{u}-\mathbf{u}_{n}\right) |,  \notag \\
& \leqslant C(||(a_{n},b_{n})||_{\mathcal{H}_{p}(\Gamma )}^{2}+||(a,b)||_{%
\mathcal{H}_{p}(\Gamma )}^{2}+\left\Vert \mathbf{u}\right\Vert
_{V}^{2})\Vert \mathbf{u}-\mathbf{u}_{n}\Vert _{2}^{2}+\frac{\nu }{2}%
\left\Vert \mathbf{u}-\mathbf{u}_{n}\right\Vert _{V}^{2}\,  \notag \\
& \leqslant C_{2}(h^{2}+\left\Vert \mathbf{u}\right\Vert _{V}^{2})\Vert 
\mathbf{u}-\mathbf{u}_{n}\Vert _{2}^{2}+\frac{\nu }{2}\left\Vert \mathbf{u}-%
\mathbf{u}_{n}\right\Vert _{V}^{2},  \label{BB111}
\end{align}%
where the function $h$ in \eqref{BB111} is given by \eqref{fh}.

Setting%
\begin{equation}
\mathbf{G}_{n}=\mathbf{G}(t,\mathbf{y}_{n}),\quad \mathbf{G}=\mathbf{G}(t,%
\mathbf{y}),\quad \mathbf{G}^{\ast }=\mathbf{G}^{\ast }(t),  \label{GS_NS1}
\end{equation}%
and using the same arguments as in the deductions of \eqref{GG1} by taking $%
C_{3}=2K$, we infer that 
\begin{align}
\Vert \mathbf{G}(t,\mathbf{y}_{n})-\mathbf{G}^{\ast }(t)\Vert _{2}^{2}&
\leqslant C_{3}\left\Vert \mathbf{u}_{n}-\mathbf{u}\right\Vert
_{2}^{2}+C\left\Vert \mathbf{a}_{n}-\mathbf{a}\right\Vert _{2}^{2}-\Vert 
\mathbf{G}-\mathbf{G}^{\ast }\Vert _{2}^{2}  \notag \\
& +2(\mathbf{G}_{n}-\mathbf{G}^{\ast },\mathbf{G}-\mathbf{G}^{\ast }).
\label{GG11}
\end{align}%
The positive constants $C_{2}$ and $C_{3}$ in (\ref{BB111}) and (\ref{GG11})
are independent of $n,$ and they may depend on the data.

Let us consider the function 
\begin{equation}
\widehat{\xi }(t)=e^{-\int_{0}^{t}\widehat{f}(s)\ ds}\quad \text{with \ }%
\widehat{f}(t)=\left[ C_{3}+\max (3C_{0},C_{2})(1+h^{2})\right] .\qquad
\label{kksi1}
\end{equation}%
Now, by applying It\^{o}'s formula to the equality (\ref{y131}), the
definition (\ref{kksi1}) of $\widehat{\xi }$, we obtain%
\begin{align}
d(\widehat{\xi }(t)||\mathbf{u}& -\mathbf{u}_{n}||_{2}^{2})+\frac{3\nu }{2}%
\widehat{\xi }(t)||\mathbf{u}-\mathbf{u}_{n}||_{V}^{2}\,dt\leqslant 2%
\widehat{\xi }(t)((\mathbf{y}_{n}\cdot \nabla )\mathbf{y}_{n}-B^{\ast }(t),%
\mathbf{u}-\mathbf{u}_{n})\,dt  \notag \\
& +2\nu \int_{\Gamma }\widehat{\xi }(t)(b-b_{n})((\mathbf{u}-\mathbf{u}%
_{n})\cdot {\bm{\tau }})\,d\mathbf{\gamma }+\widehat{\xi }(t)\Vert \mathbf{G}%
(t,\mathbf{y}_{n})-\mathbf{G}^{\ast }(t)\Vert _{2}^{2}\,dt  \notag \\
& -2\widehat{\xi }(t)\left( \mathbf{G}(t,\mathbf{y}_{n})-\mathbf{G}^{\ast
}(t),\mathbf{u}-\mathbf{u}_{n}\right) \,d{\mathcal{W}}_{t}+C\widehat{\xi }%
(t)\Vert (a_{n},b_{n})-(a,b)\Vert _{\mathcal{H}_{p}(\Gamma )}  \notag \\
& -C_{3}\widehat{\xi }(t)||\mathbf{u}-\mathbf{u}_{n}||_{2}^{2}\,dt-C_{2}%
\widehat{\xi }(t)(h^{2}+\left\Vert \mathbf{u}\right\Vert _{V}^{2})||\mathbf{u%
}-\mathbf{u}_{n}||_{2}^{2}\,dt.  \label{aps}
\end{align}%
Therefore, writing the inequality (\ref{aps}) in the integral form, taking
the expectation, and incorporating the estimates (\ref{BB111}), (\ref{GG11}%
), we infer that 
\begin{align}
\mathbb{E}(\widehat{\xi }(t)||\mathbf{u}(t)& -\mathbf{u}_{n}(t)||_{2}^{2})+%
\mathbb{E}\int_{0}^{t}\widehat{\xi }(s)\Vert \mathbf{G}-\mathbf{G}^{\ast
}\Vert _{2}^{2}ds+\nu \mathbb{\ E}\int_{0}^{t}\widehat{\xi }(s)||\mathbf{u}-%
\mathbf{u}_{n}||_{V}^{2}\,ds  \notag \\
& \leqslant 2\nu \mathbb{E}\int_{0}^{t}\int_{\Gamma }\widehat{\xi }%
(s)(b-b_{n})((\mathbf{u}-\mathbf{u}_{n})\cdot {\bm{\tau }})\,d\mathbf{\gamma 
}ds  \notag \\
& +2\mathbb{E}\int_{0}^{t}\widehat{\xi }(s)|\left( B_{1},\mathbf{u}-\mathbf{u%
}_{n}\right) |\,ds+2\mathbb{E}\int_{0}^{t}~\widehat{\xi }(s)|\left( B_{2},%
\mathbf{u}-\mathbf{u}_{n}\right) |\,ds  \notag \\
& +2\mathbb{E}\int_{0}^{t}\widehat{\xi }(s)\left( B_{3},\mathbf{u}-\mathbf{u}%
_{n}\right) \,ds+C\mathbb{E}\int_{0}^{t}\widehat{\xi }(s)\Vert
(a_{n},b_{n})-(a,b)\Vert _{\mathcal{H}_{p}(\Gamma )}^{2}\,ds  \notag \\
& +2\mathbb{E}\int_{0}^{t}\widehat{\xi }(s)(\mathbf{G}_{n}-\mathbf{G}^{\ast
},\mathbf{G}-\mathbf{G}^{\ast })\,ds  \notag \\
& =J_{0}+J_{1}+J_{2}+J_{3}+J_{4}+J_{5}\qquad \text{for \ }t\in (0,T).
\label{14:41}
\end{align}%
In the following, we show that the right-hand side of this inequality tends
to zero as $n\rightarrow \infty $. \ The H\"{o}lder inequality, (\ref%
{eq00sec12}), (\ref{cal}) and $\widehat{\xi }\leqslant \xi _{h}^{2}$ yield 
\begin{eqnarray*}
J_{0} &=&|2\nu \int_{\Gamma }\widehat{\xi }(s)(b-b_{n})((\mathbf{u}-\mathbf{u%
}_{n})\cdot {\bm{\tau }})\,d\mathbf{\gamma }ds| \\
&\leqslant &C\Vert (a_{n},b_{n})-(a,b)\Vert _{L_{2}(\Omega \times (0,T);%
\mathcal{H}_{p}(\Gamma ))}+\frac{\nu }{2}\mathbb{E}\int_{0}^{t}\widehat{\xi }%
(s)\left\Vert \mathbf{u}-\mathbf{u}_{n}\right\Vert _{V}^{2}\,ds.
\end{eqnarray*}%
Considering the estimate (\ref{lp11}) and using that $\widehat{\xi }%
\leqslant \xi _{h}^{2}$ on $(0,T)$, we deduce that

\begin{align*}
J_{1}& \leqslant \mathbb{E}\int_{0}^{T}\xi _{h}^{2}\left( ((\mathbf{u}_{n}+%
\mathbf{a}_{n})\cdot \nabla )(\mathbf{a}_{n}-\mathbf{a}),\mathbf{u}_{n}-%
\mathbf{u}\right) \leqslant C(\mathbb{E}\int_{0}^{T}\Vert \mathbf{a}_{n}-%
\mathbf{a}\Vert _{H^{1}}^{2}ds)^{1/2} \\
& \times \lbrack (\mathbb{E}\int_{0}^{T}\xi _{h}^{4}\Vert \mathbf{u}%
_{n}\Vert _{4}^{2}\Vert \mathbf{u}_{n}-\mathbf{u}\Vert _{4}^{2}ds)^{1/2}+(%
\mathbb{E}\int_{0}^{T}\xi _{h}^{4}\Vert \mathbf{a}_{n}\Vert _{C(\bar{%
\mathcal{O}})}^{2}\Vert \mathbf{u}_{n}-\mathbf{u}\Vert _{2}^{2}ds)^{1/2}] \\
& \leqslant C\Vert (a_{n},b_{n})-(a,b)\Vert _{L_{2}(\Omega \times (0,T);%
\mathcal{H}_{p}(\Gamma ))}\rightarrow 0\;\quad \text{as }n\rightarrow \infty
,
\end{align*}%
where we used the following uniform estimates with respect to the parameter $%
n$ 
\begin{align*}
\mathbb{E}& \int_{0}^{T}\xi _{h}^{4}\Vert \mathbf{u}\Vert _{4}^{2}\Vert 
\mathbf{u}_{n}-\mathbf{u}\Vert _{4}^{2}ds \\
& \leqslant (\mathbb{E}\sup_{t\in \lbrack 0,T]}\xi _{h}^{4}\Vert \mathbf{u}%
\Vert _{2}^{2}\Vert \mathbf{u}_{n}-\mathbf{u}\Vert _{2}^{2})^{\frac{1}{2}%
}\times (\mathbb{E(}\int_{0}^{T}\xi _{h}^{2}\Vert \mathbf{u}\Vert _{V}\Vert 
\mathbf{u}_{n}-\mathbf{u}\Vert _{V}ds)^{2})^{\frac{1}{2}} \\
& \leqslant C(\mathbb{E}\sup_{t\in \lbrack 0,T]}\xi _{h}^{4}\left( \Vert 
\mathbf{u}\Vert _{2}^{4}+\Vert \mathbf{u}_{n}\Vert _{2}^{4}\right) )^{\frac{1%
}{2}}\times (\mathbb{E(}\int_{0}^{T}\xi _{h}^{2}\left( \Vert \mathbf{u}\Vert
_{V}^{2}+\Vert \mathbf{u}\Vert _{V}^{2}\right) ds)^{2})^{\frac{1}{2}%
}\leqslant C,
\end{align*}%
\begin{align*}
\mathbb{E}& \int_{0}^{T}\xi _{h}^{4}\Vert \mathbf{a}_{n}\Vert _{C(\bar{%
\mathcal{O}})}^{2}\Vert \mathbf{u}_{n}-\mathbf{u}\Vert _{2}^{2}ds \\
& \leqslant C(\mathbb{E}\sup_{t\in \lbrack 0,T]}\xi _{h}^{4}\left( \Vert 
\mathbf{u}_{n}\Vert _{2}^{4}+\Vert \mathbf{u}\Vert _{2}^{4}\right) )^{\frac{1%
}{2}}\times (\mathbb{E(}\int_{0}^{T}\Vert (a_{n},b_{n})\Vert _{\mathcal{H}%
_{p}}^{2}ds)^{2})^{\frac{1}{2}}\leqslant C.
\end{align*}%
For the term $J_{2},$ using H\"{o}lder's inequality, \eqref{ksi2}, \eqref{y}
and \eqref{lp11}, we can show that 
\begin{align*}
J_{2}& \leqslant \mathbb{E}\int_{0}^{T}\xi _{h}^{2}\left( ((\mathbf{a}_{n}-%
\mathbf{a})\cdot \nabla )(\mathbf{u}+\mathbf{a}),\mathbf{u}_{n}-\mathbf{u}%
\right) \\
& \leqslant C(\mathbb{E}\int_{0}^{T}\Vert \mathbf{a}_{n}-\mathbf{a}\Vert _{C(%
\bar{\mathcal{O}})}^{2}ds)^{\frac{1}{2}} \\
& \times (\mathbb{E}\sup_{t\in \lbrack 0,T]}\xi _{h}^{2}\Vert \mathbf{u}_{n}-%
\mathbf{u}\Vert _{2}^{2}\int_{0}^{T}\left( \xi _{h}^{2}\Vert \mathbf{u}\Vert
_{V}^{2}+\Vert (a,b)\Vert _{\mathcal{H}_{p}}^{2}\right) ds)^{\frac{1}{2}} \\
& \leqslant C\Vert (a_{n},b_{n})-(a,b)\Vert _{L_{2}(\Omega \times (0,T);%
\mathcal{H}_{p}(\Gamma ))}\rightarrow 0.
\end{align*}%
Therefore, the terms $J_{1},$ $J_{2}$\ converge to zero as $n\rightarrow
\infty $.

The convergence \eqref{c11} shows that 
\begin{equation*}
\xi _{h}\left( \mathbf{u}-\mathbf{u}_{n}\right) \rightharpoonup 0\qquad 
\text{ weakly in }L_{2}(\Omega \times (0,T),V).
\end{equation*}%
The operator $\xi _{h}^{2}B_{3}=\xi _{h}^{2}\left( (\mathbf{y}\cdot \nabla )%
\mathbf{y}-B^{\ast }\right) $ belongs to $L_{2}(\Omega \times (0,T);{%
V^{\prime }})$\ by (\ref{l}), thus (\ref{c011}) implies 
\begin{equation*}
J_{3}=2\mathbb{E}\int_{0}^{T}\widehat{\xi }(s)((\mathbf{y}\cdot \nabla )%
\mathbf{y})-B^{\ast },\mathbf{u}-\mathbf{u}_{n})\,ds\rightarrow 0.
\end{equation*}%
The term%
\begin{eqnarray*}
J_{4} &=&C\mathbb{E}\int_{0}^{t}\widehat{\xi }(s)\Vert
(a_{n},b_{n})-(a,b)\Vert _{\mathcal{H}_{p}(\Gamma )}\,ds \\
&\leqslant &C\Vert (a_{n},b_{n})-(a,b)\Vert _{L_{2}(\Omega \times (0,T);%
\mathcal{H}_{p}(\Gamma ))}\rightarrow 0\;\quad \text{as }n\rightarrow \infty
.
\end{eqnarray*}%
Due to the convergence results (\ref{c011}) and (\ref{GS_NS1}), we obtain 
\begin{equation}
\xi _{h}\left( \mathbf{G}_{n}-\mathbf{G}^{\ast }\right) \rightharpoonup 
\mathbf{0}\text{\qquad\ weakly in }L_{2}(\Omega \times (0,T),H^{m}),
\label{4.19}
\end{equation}%
which implies 
\begin{equation*}
J_{5}=2\mathbb{E}\int_{0}^{T}\widehat{\xi }(s)(\mathbf{G}_{n}-\mathbf{G}%
^{\ast },\mathbf{G}-\mathbf{G}^{\ast })\,ds\rightarrow 0.
\end{equation*}

Gathering the convergence results for $J_{i}$, $i=0,\dots ,5$, and passing
to the limit in the inequality \eqref{14:41}, we deduce the following strong
convergences\ 
\begin{equation}
\lim_{n\rightarrow \infty }\mathbb{E}\left( \widehat{\xi }(t)||\mathbf{u}(t)-%
\mathbf{u}_{n}(t)||_{2}^{2}\right) =0,\text{\quad }\lim_{n\rightarrow \infty
}\mathbb{E}\int_{0}^{t}\widehat{\xi }(s)||\mathbf{u}-\mathbf{u}%
_{n}||_{V}^{2}\,ds=0  \label{gg1}
\end{equation}%
for $t\in (0,T).$ In addition, we obtain 
\begin{equation*}
\mathbb{E}\int_{0}^{t}\widehat{\xi }(s)\Vert \mathbf{G}(s,\mathbf{y})-%
\mathbf{G}^{\ast }(s)\Vert ds=0,
\end{equation*}%
then 
\begin{equation}
\mathbf{G}(t,\mathbf{y})=\mathbf{G}^{\ast }(t)\qquad \text{a. e. in }\Omega
\times (0,T).  \label{09}
\end{equation}

On the other hand, from (\ref{c011}) and (\ref{gg1}), we infer that $%
\widehat{\xi }(t)(\mathbf{y}\cdot \nabla )\mathbf{y}=\widehat{\xi }(t)%
\mathbf{B}^{\ast }(t)$ a.e. in $\Omega \times (0,T),$ that implies 
\begin{equation}
(\mathbf{y}\cdot \nabla )\mathbf{y}=\mathbf{B}^{\ast }(t)\qquad \text{a. e.
in }\Omega \times (0,T).  \label{RC2Y}
\end{equation}

Considering the identifications \eqref{09}-\eqref{RC2Y}, the equation (\ref%
{cp2}) reads 
\begin{eqnarray*}
\left( \mathbf{y}(t),\boldsymbol{\varphi }\right) -\left( \mathbf{y}_{0},%
\boldsymbol{\varphi }\right)  &=&\int_{0}^{t}\left[ -\nu \left( \mathbf{y},%
\boldsymbol{\varphi }\right) _{V}+\nu \int_{\Gamma }b(\boldsymbol{\varphi }%
\cdot {\bm{\tau }})\,d\mathbf{\gamma }-(B^{\ast }(s),\boldsymbol{\varphi })%
\right] \,ds \\
&+&\int_{0}^{t}\left( \mathbf{G}^{\ast }(s),\boldsymbol{\varphi }\right) \,d{%
\mathcal{W}}_{s},\quad P\text{-a.e. in }\Omega ,\quad t\in (0,T),
\end{eqnarray*}%
for any $\boldsymbol{\varphi }\in V.$ Therefore $\mathbf{y}$ is the solution
of the state equation, corresponding to the control pair $(a,b).$

Taking into account the lower semicontinuity of the cost functional, the
strong convergence (\ref{gg1}) and Remark \ref{R4.1}, we infer that 
\begin{equation*}
J(a,b,\mathbf{y})\leqslant \lim_{n\rightarrow \infty }J(a_{n},b_{n},\mathbf{y%
}_{n}),
\end{equation*}%
which implies 
\begin{equation*}
J(a,b,\mathbf{y})=\inf (\mathcal{P}),
\end{equation*}
hence\ the triplet $(a,b,\mathbf{y})$ is a solution to the control problem $(%
\mathcal{P})$.
\end{proof}

\bigskip

\section{Conclusion and discussion}
This work adresses an optimal control problem for the evolution of a viscous incompressible Newtonian fluid filling a two-dimensional bounded domain, under the action of  random forces modeled by a multiplicative Gaussian noise.
We prove the existence and uniqueness of the solution to the stochastic state equation and establish the existence of an optimal control.
The control is exerted at the boundary through the physical non-homogeneous  Navier-slip boundary conditions.

 Let us emphasise that
the  studies in the literature 
 \cite{C}, \cite{CC3}-\cite{CC6} turn out that the 
 non-homogeneous Navier-slip boundary conditions are  compatible with the inviscid limit transition of the viscous state, then we expect that our approach will be relevant 
to control  the evolution of turbulent flows typically  associated with  high Reynolds number (or small viscosity).

We should mention that the most results in the literature on the optimal control of 
fluid flows  are of deterministic nature. 
The control of a  stochastic system is much more involved and there are 
 few results available in the literature. We refer  the articles
   \cite{B99}, \cite{L00}  and  \cite{BT19}, \cite{BT21}, where the authors 
  solved tracking control problems  in 2D and 3D, respectively. In these works, the control variables act in the interior of the domain.
  Recently in  \cite{ZG23},
the authors studied a stochastic boundary control problem for the deterministic steady Navier-Stokes equations, where the stochastic control is imposed on the boundary by a stochastic non-homogeneous Dirichlet boundary condition.

In a forthcoming paper, we  intend to deduce the first-order necessary optimality 
conditions and analyse  the second-order sufficient conditions,  which are important for implementing numerical methods to determine the optimal boundary control.

\bigskip

\textbf{Acknowledgments.} 
{\ A substantial part of this work was developed during N.V. Chemetov's
visit to the NOVAMath Research Center. He would like to thank the NOVAMath
for the financial support (through the projects UIDB/00297/2020 and
UIDP/00297/2020) and the very good working conditions. }
The work of N.V. Chemetov was also supported by FAPESP
(Funda\c{c}\~{a}o de Amparo \`{a} Pesquisa do Estado de S\~{a}o Paulo),
project 2021/03758-8, "Mathematical problems in fluid dynamics".

The work of F. Cipriano is funded by national funds through the FCT - Funda%
\c{c}\~{a}o para a Ci\^{e}ncia e a Tecnologia, I.P., under the scope of the
projects UIDB/00297/2020 and UIDP/00297/2020 (Center for Mathematics and
Applications).

\end{document}